\documentclass[11pt, A4paper]{article}
\usepackage[margin=1in,footskip=0.25in]{geometry}
\usepackage{amsmath}
\usepackage{amsfonts}
\usepackage{amssymb}
\usepackage{amsthm,latexsym,epsfig}
\usepackage{graphicx}
\usepackage{algorithm}
\usepackage{algpseudocode}
\usepackage[colorlinks]{hyperref}
\hypersetup{
    citecolor= {blue}
}

\newtheorem{theorem}{Theorem}
\newtheorem{proposition}[theorem]{Proposition}%

\newtheorem{example}{Example}%
\usepackage{epsfig}

\usepackage{color,epsf}
\usepackage{graphicx}
\usepackage{subfig}
\usepackage{subfloat}

\newtheorem{lemma}[theorem]{Lemma}

\newcommand{\lyxdot}{.}

\newcommand{\Npperp}{N^\prime_\perp}
\newcommand{\Nppar}{N^\prime_\|}
\newcommand{\bbR}{\mathbb R}
\newcommand{\vperp}{v_\perp}
\newcommand{\vpar}{v_\|}

\begin{document}


\title{On the unconventional Hug integrator}
\author{Christophe Andrieu \footnote{School of Mathematics, University of Bristol, U.K.} and J.M. Sanz-Serna\footnote{Departamento de Matem\'aticas, Universidad Carlos III de Madrid, Avenida de la
Universidad 30, E-28911 Legan\'es (Madrid), Spain. E-mail: jmsanzserna@gmail.com}}
\maketitle

\abstract{Hug is a recently proposed iterative mapping used to design efficient updates in Markov chain Monte Carlo (MCMC) methods. Hug generates proposals that remain very close to hypersurfaces (level sets) of constant probabilty density. We analyse a generalization of Hug from hypersurfaces to manifolds of arbitrary dimensions, not necessarily arising in a sampling context.
The analysis is based on interpreting, in a nonstandard way, Hug as a consistent discretization of a system of
differential equations with a rather complicated structure. The proof of convergence of this discretization includes a number of unusual features we explore fully, in particular a supraconvergence property is established, whereby second order of convergence is attained with consistency of the first order. We uncover and discuss an unexpected  property of the solutions of the underlying dynamical system that manifest itself by the existence of Hug trajectories that fail to cover the manifold of interest.}


%

\section{Introduction}

Hug and Hop are the  building blocks of a Markov chain Monte Carlo (MCMC) method suggested in \cite{LS23} to sample from complex probability distributions \cite{youhan-serna-2014}. The method may be viewed as an alternative to the well-known HMC (Hybrid/Hamiltonian Monte Carlo) algorithm \cite{duane1987hybrid}; as HMC, it turns local information
about the target density of interest, available in the form of gradients, into an efficient mapping adapted to the
geometry of the problem at hand. The algorithm in  \cite{LS23}  has the attractive
feature that an update involves moving along contours of the probability
density of interest itself, which should be contrasted with the corresponding,
sometimes unintuitive, HMC trajectories arising from solutions of Hamilton's
equations.

The Hug update of \cite{LS23} is based on the repeated composition of
a ``bounce'' mapping (proposed by \cite{peters2012rejection} but also appearing in \cite{neal2003slice}
and exploited in \cite{bouchard2018bouncy}, \cite{vanetti2017piecewise, vanetti2019piecewise} and \cite{sherlock2021discretebouncyparticlesampler,sherlock2022discrete}) in both the continuous and discrete time Markov chain contexts). The bouncing mechanism makes it possible for the proposals generated by Hug to remain very close to (hug) a manifold of interest, which is remarkable since the algorithm is fully explicit, i.e.\ it does not require solving the algebraic equations defining the manifold, in contrast with exact methods e.g. \cite{zappa-holmes}.  This combination of favourable properties makes Hug potentially relevant to problems beyond sampling and calls for a better understanding of its properties.

The aim of this paper is to provide a theoretical analysis of Hug. We clarify from the outset that comparing the practical performance of Hug with alternative techniques, both in sampling and nonsampling contexts,
 is beyond the scope of our study.  Rather our aim is to develop tools that may be useful to understand the properties of algorithms based on Hug. More specifically we establish that Hug may be viewed as a consistent discretization of an underpinning dynamical system, which we identify.  In the HMC algorithm context, where the Leapfrog integrator arises from Hamilton's equations, consistency of the discretisation plays an essential role in developing the scaling limits arguments that have led to the practical implementational rules on how to optimally tune the parameters of the algorithm \cite{beskos-pillai-2013}.

The contents and contributions of this manuscript are as follows. Section~\ref{sec:Hug} focuses on discrete-time algorithms. In Subsection~\ref{subsec:originalHug} we
present the Hug update in its original MCMC context where the aim
is to maintain a single constraint while updating the state of the Markov chain.
In Subsection~\ref{subsec:generalHug} we present a generalisation of Hug's discrete
dynamics \cite{C24} suitable for the exploration of manifolds resulting from
multiple constraints; the need to define dynamical systems that remain on a manifold arises in contexts beyond the sampling context (see Example~\ref{ex:splitting})  \cite[Chapter 7]{Leimkuhler_Reich_2005}.  A number of properties of the algorithm are presented in Subsection~\ref{subsec:properties}. In Section~\ref{sec:ODE}
we identify the  (rather unusual) underpinning continuous-time dynamical system and
explain the  nonstandard way in which Hug provides a convergent discretization of that system. In Subsection~\ref{subsec:differential-projector}
we review results on linear projector differentials. Then, in Subsection~\ref{subsec:ODE-the-equations} we present the continuous-time dynamical system  and establish a number of properties that mirror those of the discrete time algorithm studied previously in Subsection~\ref{subsec:properties}.
The convergence of
Hug as a discretization of the continuous-time dynamical system is elucidated in Subsection~\ref{subsec:relation-ODE-integrator}. The analysis presents a number of nonstandard features; of particular interest from the numerical analysis point of view is that a  \emph{supraconvergence} phenomenon \cite{GS91} takes place, whereby second-order convergence is obtained with first order consistency.
In Subsection~\ref{subsec:worked-out-example} we fully analyse a model system, showing that, unexpectedly, the trajectories generated by both the differential equations and Hug may fold back on themselves, thus precluding, in some conditions, good exploration of the manifold of interest. This unwelcome phenomenon depends on the initial condition and we obtain
a full characterisation of the initial states for which the phenomenon occurs. In Section~\ref{sec:Numerical-experiments-and} we present simple numerical experiments
suggesting that the phenomenon analysed in the model problem persists in more complex scenarios. We also comment briefly on how to suppress or mitigate the folding-back phenomenon.
Technical proofs are given in Section~\ref{sec:proofs}; these make extensive use of the results  by Golub and Pereyra \cite{GP} on differentiation of projections.

\section{Hug algorithm} \label{sec:Hug}
\subsection{The algorithm} \label{subsec:originalHug}
In this subsection we describe the Hug algorithm suggested in \cite{LS23}. In the description we will make use of some properties that will be proved in Proposition~\ref{prop:alg} below.

Hug is a Markov Chain Monte Carlo method to obtain samples from a target unnormalized probability distribution \(\pi\) in \(\bbR^n\), \(n>1\), assumed
to have a density \(\exp(\ell(x))\) with respect to the Lebesgue measure. Hug introduces an auxiliary variable \(v\in\bbR^n\) with a density \(q(v|x)\) such that
\begin{equation}\label{eq:qreversible}
q(v|x)=q(-v|x)
\end{equation}
and generates a Markov Chain with invariant distribution \(\pi(x)q(v|x)\). The marginal on \(x\) of this Markov chain has therefore \(\pi(x)\) as an invariant distribution;  one step of this marginal chain is described in Algorithm~\ref{alg:hug}.

A salient feature of the algorithm is the use of velocity reflections. For \(x\in\bbR^n\) with \(\nabla\ell(x)\neq 0\), the \(n\times n\) matrix \(R(x)\) used in the reflections is
\begin{equation}\label{eq:R}
R(x) = I-\frac{2}{\| \nabla \ell(x)\|^2} \nabla \ell(x) \nabla \ell(x)^{\top},
\end{equation}
so that, if \(w\) is a vector, \(R(x)w\) is the result of reflecting \(w\) on the vector subspace of \(\bbR^n\) tangent at \(x\) to the level hypersurface \(\{y:\ell(y)=\ell(x)\}\). Here and elsewhere \(\|\cdot\|\) denotes the standard Euclidean norm.

\begin{figure}[t]
	\centering
	\includegraphics[width=0.6\textwidth]{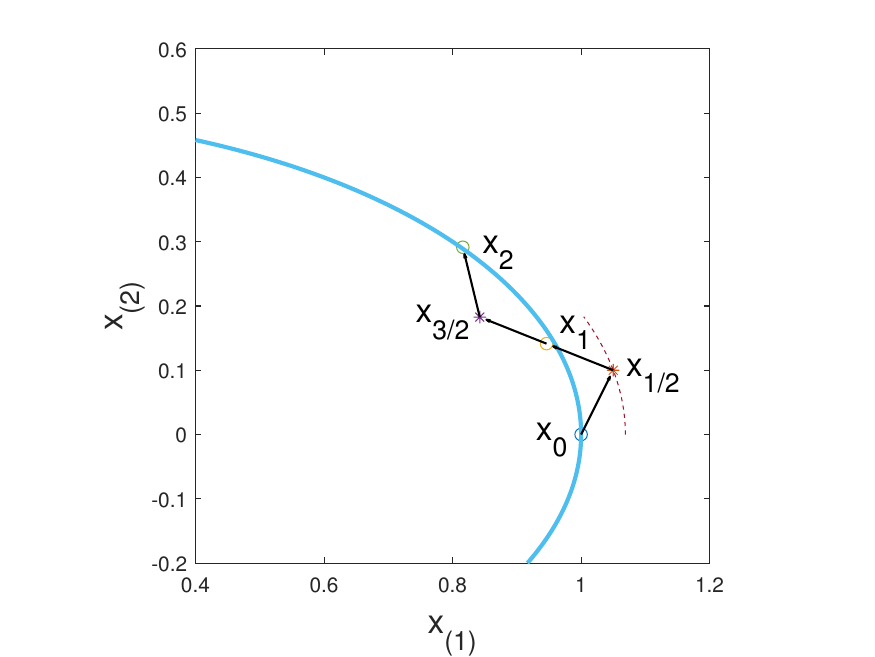}
	\caption{Computation of the iterate \(x_2\) for a bivariate Gaussian. The plot depicts the plane of the variable \(x=(x_{(1)},x_{(2)})\in\bbR^2\). The solid curve is the contour  \(\{y :\ell(y) = \ell(x_0)\}\).  The dashed curve is (part of) the contour  \(\{y :\ell(y) = \ell(x_{1/2})\}\). The segments \([x_0,x_{1/2}]\), \([x_{1/2},x_{1}]\), \([x_1,x_{3/2}]\), \([x_{3/2},x_{2}]\)
share a common Euclidean length \((\delta/2)\|v_0\|\).}
	\label{fig:hug}
\end{figure}

Figure~\ref{fig:hug} illustrates the computation of the iterate \(x_2\) when \(\delta=0.1\), \(x_0=[1,0]^{\top}\), \(v_0=[1,2]^{\top}\) and the target is a bivariate Gaussian distribution with
\(\ell = -x_{(1)}^2-4x_{(2)}^2\) (\(x_{(1)}\) and \(x_{(2)}\) denote the scalar components of \(x\)).  The figure suggests that  the iterates \(x_k\) may be expected  to be close to (or \lq\lq hug\rq\rq)  the level set \(\{y:\ell(y)=\ell(x_0)\}\); this explains the name of the algorithm. The  halfway points \(x_{k+1/2}\) will not be too close to \(\{y:\ell(y)=\ell(x_0)\}\), but this does not matter as they are just intermediate auxiliary values (like the internal stages of a Runge-Kutta integrator).

\begin{algorithm}[t]
	\caption{Hug} \label{alg:hug}
\begin{algorithmic}[b]
		\Require Number of timesteps \(K\); stepsize \(\delta>0\); current state \(x\) of the Markov chain; \(v\)-marginal density  \(q(\cdot|x)\).

         Initialize \(k = 0\), \(x_0\leftarrow x\), draw velocity \(v_0\sim q(\cdot|,x_0)\).

		  Timestepping: \For{\(k=0,\dots,K-1,\)}

            Move to \(x_{k+1/2} = x_k +(\delta/2) v_k\).

            Reflect: \(v_{k+1} = R(x_{k+1/2})v_k\).

            Move to \(x_{k+1} = x_{k+1/2} +(\delta/2) v_{k+1}\).

          \EndFor
           		
		Compute \(\log(r) = \ell(x_K)-\ell(x_0)+\log q(v_K|x_K)-\log q(v_0|x_0)\).

    With a probability \(\alpha = 1\wedge r\), \(x\leftarrow x_K\); else \(x\leftarrow x\).
\end{algorithmic}
\end{algorithm}

The typical choice of \(q(v|x)\) satisfying \eqref{eq:qreversible} is given by an \(x\)-independent isotropic normal \(v\sim {\mathcal N}(0,\sigma^2 I)\). Since the reflection matrix \eqref{eq:R} is orthogonal, all iterates \(v_0\), \dots, \(v_K\) share a common Euclidean length, and therefore, for this simple choice of \(q\), in the expression for the acceptance probability \(r\) the term \[\log q(v_K|x_K)-\log q(v_0|x_0)= -\sigma^{-2} \big(\|v_K\|^2-\|v_0\|^2\big)\] vanishes.
Furthermore, due to the \lq\lq hugging\rq\rq\ property,  \(\ell(x_K)-\ell(x_0)\) will typically be small so that the value of the acceptance probability \(r\) will be close to \(1\).
At the same time, it could be expected that, once \(\delta\) has been chosen, by taking the number \(K\) of timesteps sufficiently high, \(x_K\) will be far away from \(x_0\). In conclusion, Hug may offer the possibility of generating proposals that are away from the current state of the chain and, at the same time, have high probability of being accepted.

The \(x\)-moves in Hug will not change much the value of \(\ell\) and in order to efficiently explore the whole state space,
the reference \cite{LS23} suggested to interleave steps of Hug with steps of a second Markov kernel, called Hop, that causes the state of the Markov chain to jump between different level sets   of \(\ell\).

\subsection{A generalization of the Hug timestepping mechanism} \label{subsec:generalHug}
In what remains of the paper, we will consider that in the Hug time-stepping formulas (\(x_0\in\bbR^n\), \(v_0\in\bbR^n\), \(n\geq 2\))
\begin{eqnarray}
x_{k+1/2} &=& x_k +(\delta/2) v_k\label{eq:alg1}\\
v_{k+1} &=& R(x_{k+1/2})v_k,\label{eq:alg2}\\
x_{k+1} &=& x_{k+1/2} +(\delta/2) v_{k+1},\label{eq:alg3}
\end{eqnarray}
the reflection in \eqref{eq:alg2} is taken with respect to  the vector subspace of \(\bbR^n\) tangent at \(x_{k+1/2}\) to the level manifold \(\{y:f(y)=f(x_{k+1/2})\}\) of a smooth function \(f:\bbR^n\rightarrow\bbR^m\), \(0<m<n\). It is expected that now the iterates
\(x_k\) will hug the level set \(\{y:f(y)=f(x_0)\}\). To shorten the notation we shall hereafter set \({\mathcal M}(\eta) = \{y:f(y)=\eta\}\)  \((\eta\in\bbR^m\)).
In the particular case of Algorithm~\ref{alg:hug}, \(m=1\) and
\(f=\ell\). The present generalization has been considered in the thesis \cite{C24}, where it is used to build algorithms to sample from \emph{filamentary distributions}, i.e.\ distributions in \(\bbR^n\) supported in the neighbourhood of a set \({\mathcal M}(\eta)\) for some fixed  \(\eta\in\bbR^m\). An application of the \(m>1\) case to  a problem not related to sampling  is outlined in Example \ref{ex:splitting}.

We now present the expression for \(R(x)\) that replaces \eqref{eq:R} in the generalized scenario. For each \(x\in\bbR^n\), \(J(x)\) will denote the \(m\times n\) Jacobian matrix of \(f\) evaluated at \(x\). In order not to clutter the exposition with unwelcome details, we will assume throughout that, for each \(x\), \(J(x)\) has full rank. (However this hypothesis may easily be weakened: for instance one might demand that \(J(x)\) be full rank except at some isolated points or, more generally, that \(J(x)\) be full rank for \(x\) in an open subdomain of \(\bbR^n\).) The implicit function theorem shows that, for each \(x\), \({\mathcal M }(f(x))\) is  a smooth \(n-m\) submanifold of \(\bbR^n\). The vector space \emph{tangent} at \(x \) to this manifold  is given by
\begin{equation}\label{eq:orthogonal}
 {\mathcal T}(x) = \{v\in \bbR^n: J(x)v = 0\}.
\end{equation}
We shall denote by \({\mathcal N}(x)\subset \bbR^n\) the corresponding \emph{normal} space, i.e.\ the \(m\)-dimen\-sional vector subspace of \(\bbR^n\)   orthogonal to \({\mathcal T}(x)\).
For each \(x\), \(J(x)J(x)^{\top}\) is invertible and
\begin{equation}\label{eq:J+}
J(x)^+ = J(x)^{\top}(J(x)J(x)^{\top})^{-1}
\end{equation}
is the  Moore-Penrose inverse of \(J(x)\) (size \(n\times m\)). The \(n\times n\) symmetric matrix
\begin{equation}\label{eq:N}
N(x) = J(x)^+J(x) = J(x)^{\top} (J(x)J(x)^{\top})^{-1}J(x)
\end{equation}
represents the orthogonal projection of \(\bbR^n\) onto  \({\mathcal N}(x)\) as $N^2(x)=N(x)$ and $J(x)v=0$ implies $N(x)v=0$.  Then we define
\begin{equation}\label{eq:T}
T(x)=I-N(x)
\end{equation}
 the complementary projection onto \({\mathcal T}(x)\), as $T^2(x)=T(x)$ and $N(x)T(x)=0$. In what follows we will use repeatedly the fact that, for any given \(x\), each vector \(w\) may be decomposed as \(T(x)w+N(x)w\) into a tangential component and a normal component.

 After these preparations, we see that the expression for \(R(x)\) to be used in \eqref{eq:alg2} is \(R(x)=I-2N(x)\).
For given \(x\) and \(w\), there are different ways of computing \(N(x)w\). One that suggests itself uses the decomposition \(J(x)^{\top}=QR\), where \(Q\) is \(n\times m\) with orthonormal columns and \(R\) is \(m\times m\) upper triangular with positive diagonal entries. The columns of \(Q\) are an orthonormal basis of \({\mathcal N}(x)\) and therefore \(N(x)w = QQ^{\top}w\). In the particular case \(m=1\), \(f=\ell\), this amounts to finding
the normalized gradient \((1/\|\nabla \ell(x)\|)\nabla \ell(x)\) as in \eqref{eq:R}.
\subsection{Properties of the timestepping mechanism}
\label{subsec:properties}

Proposition~\ref{prop:alg} below summarizes some properties of \eqref{eq:alg1}--\eqref{eq:alg3}. It requires some notation.

We denote by \(\Psi_\delta:\bbR^n\times \bbR^n\rightarrow \bbR^n\times \bbR^n\) the map that carries out a single step of \eqref{eq:alg1}--\eqref{eq:alg3}, i.e.\ \(\Psi_\delta(x,v)= (x^\prime,v^\prime)\) if \(v^\prime = R(x+(\delta/2)v)v\), \(x^\prime = x+(\delta/2)v+(\delta/2)v^\prime\).

The symbol \(H(x)\) will denote the second derivative of \(f\) evaluated at \(x\); this is the symmetric bilinear map \(\bbR^n\times\bbR^n\rightarrow \bbR^m\) such that
\[
f(x+w)= f(x) + J(x)w+\frac{1}{2} H(x)[w,w]+ o(\|w\|^2), \qquad w\rightarrow 0.
\]
Concretely with $f:=(f_1,f_2,\ldots,f_m)$ then $\big[H(x)[u,w]\big]_{i} = u^{\top} \nabla^2 f_i (x) w$ for $i=1,\ldots,n$ where $\nabla^2 f_i (x)$ is the Hessian of $f_i$ at $x$.
For a symmetric bilinear map \(S:\bbR^n\times\bbR^n\rightarrow \bbR^m\),  the operator norm  is defined as
\[
\|S\| = \sup \{\|S[w_1,w_2]\|: \|w_1\|\leq 1,  \|w_2\|\leq 1\}.
\]
With this notation in place we have: one can establish the following properties. Properties 1-5 are fairly straightforward to establish but Property 6 requires more care.
\begin{proposition}\label{prop:alg}
The following properties hold:
\begin{enumerate}
\item The map \(\Psi_\delta\) is volume preserving.
\item The map \(\Psi_\delta\) is time-reversible, i.e.\ for each \(x\) and \(v\), \(\Psi_\delta(x,v) = (x^\prime,v^\prime)\) if and only if
\(\Psi_\delta(x^\prime,-v^\prime) = (x,-v)\).
\end{enumerate}
Furthemore, let \eqref{eq:alg1}--\eqref{eq:alg3} hold for \(k=0,\dots, K-1\). Then:
\begin{enumerate}
\item[3.]
For each \(k= 0,\dots,K-1\), \[T(x_{k+1/2})v_{k+1}= T(x_{k+1/2})v_{k},\quad N(x_{k+1/2})v_{k+1}= -N(x_{k+1/2})v_{k}\]
and \[ v_k+v_{k+1} = 2T(x_{k+1/2}) v_k.\]
\item[4.] The vectors \(v_0\),\dots, \(v_K\) share a common Euclidean length \(\|v_0\|\).
\item[5.] The vectors \(x_{k+1/2}-x_k\), \(k=0,\dots, K-1\), and  \(x_{k+1}-x_{k+1/2}\), \(k=0, \dots,K-1\), share a common Euclidean length
\((\delta/2)\|v_0\|\). The vectors \(x_{k+3/2}-x_{k+1/2}\), \(k = 0, \dots, K-2\) share a common Euclidean length \(\delta\|v_0\|\).
\item[6.]   If, as \(x\) ranges in \(\bbR^n\),\(\|H(x)\|\) is upper bounded by \(\beta\) and \(H(x)\) is \(\gamma\)-Lipschitz continuous, i.e.\ for each \(x,y\in\bbR^n\),
\[
\|H(x)-H(y)\| \leq \gamma \|x-y\|,
\]
we have
\[
\big\|f(x_K)-f(x_0)\big\| \leq \frac{\delta^2}{12} \|v_0\|^2 \big(3\beta+ \gamma (K-1)\delta \|v_0\|\big).
\]
\item[7.] If \(H(x)\) is of the form \ \(H(x)[w,w] =  \|w\|^2\nu\), \(\nu\in\bbR^m\), then \(f(x_k)=f(x_0)\), \(k= 0, \dots,K\).
\end{enumerate}
\end{proposition}
To prove Property 1., note that \((x,v)\mapsto \Psi_\delta(x,v) = (x^\prime,v^\prime)\) is the composition of three volume preserving maps,
\(
(x,v)\mapsto (x+(\delta/2)v,v)
\), \((x,v)\mapsto (x,R(x)v), (x,v) \mapsto (x+(\delta/2)v,v)\). Checking Property 2.\ is a trivial computation. Property 3.\ just restates \eqref{eq:alg2}.
 Property 4.\ follows from Property 3. Property 5.\ is a consequence of \eqref{eq:alg1}, \eqref{eq:alg3} and Property
 4. The proofs of the last two properties are detailed in Subsection~\ref{subsec:proof-prop-1}

Even though \(\Psi_\delta\) is volume-preserving (Property 1.), a straightforward computation shows that it  is not symplectic (the condition in e.g.\ \cite[Definition VI, 2.1]{geom} is not satisfied).

For the particular case \(f=\ell\), Properties 1.\ and 2.\ in tandem with standard results on Markov Chain Monte Carlo methods ensure that  Algorithm~\ref{alg:hug} preserves the distribution with density \(\propto \exp(\ell(x))\).

In Property 6 note that when \(\beta=0\)  the manifold \({\mathcal M}(f(x_0))\) is flat (an affine subspace); thus \(\beta\) is a measure of the curvature of that manifold. The bound in Property 6 shows that, for the deviation \(\| f(x_1)-f(x_0)\|\) to be small, \(\delta\) should be chosen in such a way that \(\delta \|v_0\|=2 \|x_{1/2}-x_0\|\) is small relative to the curvature of the manifold: a very intuitive result (see Fig.~\ref{fig:hug}). It is conceivable that by replacing the first-order, Euler-like substeps \eqref{eq:alg1}--\eqref{eq:alg3} by more sophisticated approximations, as suggested by a reviewer, it would be possible to improve the \(\mathcal{O}(\delta^2)\) bound in Property 6 to become \(\mathcal{O}(\delta^k)\), \(k>2\). 

Property 7.\ holds in particular if in the context of Algorithm~\ref{alg:hug} the target is an isotropic Gaussian. In that situation, if \(q(v|x)\) is also an isotropic Gaussian, Algorithm~\ref{alg:hug}  accepts all the proposals.

It would be possible to formulate a localized version of Property 6.,  where rather than demanding that \(H\) is  bounded by $\beta$ and $\gamma$-Lipschitz globally, those hypotheses are only demanded in a domain that contains the iterates \(x_k\).
Property 6.\ appears, in a slightly different form, in \cite{LS23} (for the case \(f=\ell\)) and  in \cite{C24} (for the general case). It shows that indeed \eqref{eq:alg1}--\eqref{eq:alg3} produces points \(x_k\) that, for \(\delta\) small, remain close to the manifold \({\mathcal M}(f(x_0))\). It is remarkable that for an \emph{explicit} timestepping algorithm where no implicit equations are solved, the deviation \(f(x_K)-f(x_0)\) grows at most \emph{linearly} with \(K\delta\), the accumulated time after taking \(K\) steps of length \(\delta\) (typically integration errors accumulate exponentially). For this reason Hug could be used in any situation not related to sampling where it is of interest to obtain points close to the manifold  \({\mathcal M}(\eta)\); an example, out of many possible, follows.

\begin{example}\label{ex:splitting}The paper \cite{blanes2013} is devoted to the construction of high accuracy splitting methods tailored to the integration of the Solar System. The most
accurate method in \cite{blanes2013}, labelled AB1064, is found by considering a family of integrators depending on ten parameters \(a_j\), \(
b_j\), \(j = 1,\dots,5\), and imposing  nine relations on the parameters, \(f(x) =0\), \(f:\bbR^{10}\rightarrow \bbR^9\), \(x = (a_1,\dots,b_5)\), that annihilate relevant terms in the  expansion of the local error in powers of the stepsize \(\delta\) and a small parameter \(\epsilon\) of physical origin. In this way, attention is restricted to a manifold of codimension one (curve) \(\mathcal M\) in the \(\bbR^{10}\) space of the parameters. AB1064 is then defined by the point on \(\mathcal M\) of minimum squared Euclidean length \(L^2=|a_1|^2+\cdots+|b_5|^2\). Finding the optimal parameter value is far from trivial
and the authors of \cite{blanes2013} resort to a sophisticated technique. As an alternative, since an initial point \(x_0\) on \(\mathcal M\)  is readily available, we may use Hug to get
 an approximate description of \(\mathcal M\) by points \(x_0\), \dots, \(x_K\). The point \(x_k\) of minimum \(L^2\), while not exactly on \(\mathcal M\), provides a good starting location to solve,  by means of
 a Newton iteration, the \(10\times 10\) system
\(f=0\), \(\nabla L^2=0\) satisfied by the coefficients we seek.
\end{example}

It is perhaps of interest to point out that applications of Hug similar to the one we have just sketched only require the knowledge of the mapping \(f\) that describes the constraints. This should be compared with situations where one knows a system of ODEs in \(\bbR^n\) whose solutions stay on a manifold \(f=0\) and the system is integrated with an (invariably implicit) algorithm devised to generate points on the manifold \cite[IV]{geom}. A well-known example is provided by the integrators SHAKE and RATTLE used in molecular dynamics \cite[VII.1.4]{geom}.

\section{The ODE} \label{sec:ODE}
Formulas \eqref{eq:alg1}--\eqref{eq:alg3} are clearly reminiscent of a discretization of a system of ODEs. In this section we shall prove that, in fact, \eqref{eq:alg1}--\eqref{eq:alg3} provide an integrator for a system that we will identify. The original reference \cite{LS23} did not discuss the connection between Hug and differential equations.

\subsection{The derivative \(N^\prime(x)\)} \label{subsec:differential-projector}

The expression for the system of ODEs approximated by \eqref{eq:alg1}--\eqref{eq:alg3} includes the derivative
\(N^\prime(x)\) of the projector \(N(x)\) defined in \eqref{eq:N}. This derivative is studied next.

For fixed \(x\in\bbR^n\), \(N^\prime(x)\) maps \emph{linearly} each vector \(w\in\bbR^n\) into an \(n\times n\) matrix, that we will denote by \(N^\prime(x)[w]\), in such a way that \(N(x+w) = N(x)+N^\prime(x)[w]+o(\|w\|)\) as \(w\rightarrow 0\). If \(z\in\bbR^n\), the notation \(N^\prime(x)[w]z\) will be used, as expected, to refer to the \(n\)-vector obtained by multiplying the matrix \(N^\prime(x)[w]\) and \(z\). It is important to note that for each fixed \(x\), \(N^\prime(x)[w]z\) depends \emph{bilinearly} on \(w\) and \(z\).

According to \cite{GP}, for fixed \(x\) and \(w\),
\begin{equation}\label{eq:Nprime}
N^\prime(x)[w] = \Npperp(x)[w]+\Nppar(x)[w],
\end{equation}
where, the \(n\times n\) matrices \(\Npperp(x)[w]\), \(\Nppar(x)[w]\) are given by
\begin{eqnarray}\label{eq:nprimeone}
\Npperp(x)[w] &=& J(x)^+ H(x)[w,\cdot] T(x),\\\label{eq:nprimetwo}
\Nppar(x)[w]   &= & \big(N^\prime_\perp(x)[w]\big)^{\top} = T(x)\big(H(x)[w,\cdot] \big)^{\top} \big(J(x)^+\big)^{\top}.
\end{eqnarray}
Here, \(J(x)^+\) and \(T(x)\) were defined in \eqref{eq:J+} and \eqref{eq:T} respectively and \(H(x)[w,\cdot]\) is the \(m\times n\) matrix that corresponds to the linear operator \(z\in\bbR^n\mapsto H(x)[w,z]\in\bbR^m\) resulting from freezing at \(w\) the first argument of the bilinear operator \(H(x)[\cdot,\cdot]\). In this way, the right hand-side of \eqref{eq:nprimeone} is the product of an \(n\times m\) matrix, an \(m\times n\) matrix and \(n\times n\) matrix.

The following result, that explains the subscripts \(\perp\) and \(\|\), gives  properties of the matrices \(\Npperp(x)[w]\), \(\Nppar(x)[w]\) that we will use repeatedly in the proofs of the results.

\begin{lemma}
\label{lemma:nprime}For each \(x\) and \(w\):
\begin{enumerate}
\item The image subspace of \(\Npperp(x)[w]\) is contained in \({\mathcal N}(x)\).
\item The kernel of \(\Npperp(x)[w]\) contains \({\mathcal N}(x)\).
\item The image subspace of \(\Nppar(x)[w]\) is contained in \({\mathcal T}(x)\).
\item The kernel of \(\Nppar(x)[w]\) contains \({\mathcal T}(x)\).
\end{enumerate}
\end{lemma}

It is useful to present an illustration:

\begin{example}\label{ex:gaussian} For the case \(n=2\), \(m=1\), \(f(x) = -ax_{(1)}^2-bx_{(2)}^2\), \(a,b>0\) (corresponding to a bivariate Gaussian in Algorithm~\ref{alg:hug}), if \(w=[w_{(1)},w_{(2)}]^{\top}\), one finds
\[
\Npperp(x)[w] =ab \frac{w_{(2)}x_{(1)}-w_{(1)}x_{(2)}}{(a^2x_{(1)}^2+b^2x_{(2)}^2)^2}
\begin{bmatrix}ax_{(1)}\\ bx_{(2)}
\end{bmatrix}
\begin{bmatrix}-bx_{(2)}& ax_{(1)}
\end{bmatrix}.
\]
This is a rank one, \(2\times 2\) matrix with image spanned by the vector
\([ax_{(1)}, bx_{(2)}]^{\top}\), whose direction coincides with that of \(\nabla f(x)\) (normal to the contour of \(f\) that contains \(x\)). The kernel is also spanned by \(\nabla f\). Note that the dependence on \(w\) is linear.

The matrix \(\Nppar(x)[w] \) is obtained by transposition. Its image and kernel are orthogonal to \(\nabla f(x)\).
\end{example}

\subsection{The ODE being approximated} \label{subsec:ODE-the-equations}

Algorithm  \eqref{eq:alg1}--\eqref{eq:alg3} will be shown  to approximate, in a sense to be specified later, the following system   of ODEs in \(\bbR^n\times \bbR^n\):\footnote{The ---rather unexpected--- expression for this system was derived by assuming the ansatz \eqref{eq:definition}, where \(x\) and \(v= \vpar+\vperp\) are solutions of a system of ODEs to be determined, substituting the ansatz in \eqref{eq:sigma}--\eqref{eq:tau} and imposing \(\sigma_k= \mathcal{O}(\delta^2)\) and \(\tau_k= \mathcal{O}(\delta^2)\) (cf. the proof of
 Theorem~\ref{th:consistency}). }
\begin{eqnarray}
\label{eq:odebis1}
\frac{d}{dt} x &=& T(x)v,\\
\label{eq:odebis2}
\frac{d}{dt} v &=&\Big(\Nppar(x)\big[(T(x)-N(x))v\big]-\Npperp(x)\big[(T(x)-N(x))v\big]\Big)v.
\end{eqnarray}

In order to shorten the notation, we will  introduce the symbols
\[\vpar = T(x)v,\qquad \vperp = N(x)v
\]
for the tangential and normal components of the velocity \(v\).
With this notation and taking into account the properties of the kernels of the matrices \(\Nppar[\vpar-\vperp]\) and \(\Npperp[\vpar-\vperp]\) in Lemma~\ref{lemma:nprime}, \eqref{eq:odebis2} may be written as
\begin{equation}\label{eq:odebis3}
\frac{d}{dt} v = \Nppar(x)[\vpar-\vperp]\vperp-\Npperp(x)[\vpar-\vperp]\vpar.
\end{equation}

Some properties of \eqref{eq:odebis1}--\eqref{eq:odebis2} are contained in the following result. Clearly 1., 2.\ and 3.\
here are exact counterparts of Items 1., 2. and 4. in Proposition~\ref{prop:alg}. In addition Item 6.\ in that proposition is an approximate version of Item 4. below.
\begin{theorem}\label{th:odebis}
The system \eqref{eq:odebis1}--\eqref{eq:odebis2} has the following properties:
\begin{enumerate}
\item It preserves volume in \(\bbR^n\times\bbR^n\).
\item It is time-reversible, i.e.\ it remains invariant after changing \(t\) into \(-t\) and \(v\) into \(-v\).
\item It conserves the Euclidean length of \(v(t)\): \((d/dt) \|v(t)\|^2=0\).
\item  \(x(t)\) remains on \({\mathcal M}(f(x(0))\): \((d/dt)f(x(t))=0\)
\end{enumerate}
\end{theorem}

The next result gives  expressions for \((d/d)\vpar\) and \((d/dt)\vperp\). Note that  Lemma~\ref{lemma:nprime} implies that in \eqref{eq:ode2} the terms \(- \Nppar(x)[\vperp]\vperp\)  and \(-\Npperp(x)[\vpar]\vpar\) are respectively the tangential and normal components of \((d/dt) \vpar\). Similarly the tangential and normal components of \((d/dt) \vperp\) are
\(\Nppar(x)[\vpar]\vperp\) and \(\Npperp(x)[\vperp]\vpar\).
\begin{proposition}
\label{prop:components}If \((x(t),v(t))\) is a solution of \eqref{eq:odebis1}--\eqref{eq:odebis2}, then:
\begin{eqnarray}
\label{eq:ode2}
\frac{d}{dt} \vpar &=&- \Nppar(x)[\vperp]\vperp-\Npperp(x)[\vpar]\vpar,\\
\label{eq:ode3}
\frac{d}{dt} \vperp &=& \Nppar(x)[\vpar]\vperp+\Npperp(x)[\vperp]\vpar.
\end{eqnarray}
\end{proposition}
\subsection{Relating the Algorithm and the ODE} \label{subsec:relation-ODE-integrator}

Let us now study the relation between the Algorithm \eqref{eq:alg1}--\eqref{eq:alg3} and the system \eqref{eq:odebis1}--\eqref{eq:odebis2}.
Before we do that, it is convenient to rewrite the algorithm after elimination of the auxiliary halfway points \(x_{k+1/2}\):
\begin{eqnarray}\label{eq:algbis1}
x_{k+1} &=& x_k+\delta T(x_k+(\delta/2)v_k)v_k,\\\label{eq:algbis2}
v_{k+1} &= & \big(I-2N(x_k+(\delta/2)v_k)\big)v_k.
\end{eqnarray}
The situation is much complicated by the way the algorithm
treats the normal component of the velocities \(v_k\). According to \eqref{eq:algbis2}, in the limit \(\delta\rightarrow 0\) with fixed \(x_0\), \(N(x_0)v_1=-N(x_0)v_0\); therefore the \(v_k\)'s cannot be seen as approximations to the values of a differentiable function at the step points \(k\delta\).
We will show presently that the values \(x_k\), \(v_k\) generated by the algorithm approximate the values, \(k=0,\dots,K\)
\begin{equation}\label{eq:definition}
X_k = x(k\delta),\qquad V_k = \vpar(k\delta)+(-1)^k\vperp(k\delta),
\end{equation}
where \((x(t),v(t))\) is a suitable solution of \eqref{eq:odebis1}--\eqref{eq:odebis2}. The factor \((-1)^k\)
makes it possible for the \(V_k\)'s to mimic the reflections in the normal components of the \(v_k\)'s.

The first step in the analysis is to show \emph{consistency}, i.e.\ that
when \(X_k\), \(V_k\) are substituted into the equations \eqref{eq:algbis1}--\eqref{eq:algbis2}, they originate small residuals
\begin{align}\label{eq:sigma}
\sigma_{k+1} &:= X_{k+1}- X_k-\delta T(X_k+(\delta/2)V_k)V_k,\\
\label{eq:tau}
\tau_{k+1} &:= V_{k+1} - \big(I-2N(X_k+(\delta/2)V_k)\big)V_k.
\end{align}

As with any other explicit timestepping algorithm, these residuals may be seen as truncation errors (also referred to as local errors), i.e.\ as the error of the algorithm at \(t=(k+1)\delta\) if started from the exact values at \(t=k\delta\). More precisely,
if for some \(k\), \(x_k\) and \(v_k\) happened to coincide with \(x(k\delta)\) and \(\vpar(k\delta)+(-1)^k\vperp(k\delta)\) respectively, then the algorithm would deliver, after one step, values \[x_{k+1}=X_{k+1}-\sigma_{k+1}, \qquad v_{k+1}=V_{k+1}-\tau_{k+1}.\]

The following result shows consistency of order one:
\begin{theorem}\label{th:consistency}\emph{(First-order consistency.)}
As \(\delta\rightarrow 0\) and \(k\rightarrow \infty\) with \(k\delta \rightarrow t\), the truncation errors \(\sigma_k\),\(\tau_k\) are
\(\mathcal{O}(\delta^2)\).
\end{theorem}

Consistency is the key step to proving the next convergence result. It is remarkable that, with \emph{first order} consistency, the algorithm is convergent of the \emph{second order}. This phenomenon, known as \emph{supraconvergence} is discussed  in e.g.\ \cite{GS91}, where additional references are also supplied. Supraconvergence happens here because there is substantial cancellation between  the truncation errors at consecutive steps. A numerical illutration of this cancellation is provided
in Table~\ref{table:cancellation}, where \(f\) is as in Figure \ref{fig:hug}, \(x_0 = [\cos(1), (1/2)\sin(1)]^{\top}\), \(v_0= [0,1]^{\top}\). We see that the errors after a single step (i.e.\ the truncation errors) decrease as \(\mathcal{O}(\delta^2)\) in agreement with Theorem~\ref{th:consistency}. However when a second step is taken with the same value of \(\delta\), the error decreases substantially. Errors after two steps exhibit a \(\mathcal{O}(\delta^3)\) behaviour. These issues are discussed in detail in Section~\ref{sec:proofs}. For the time being we emphasize that this cancellation of the local error in consecutive steps is \emph{not} the zig-zagging of the midway approximation \(x_{k+1/2}\) that may be seen in Figure~\ref{fig:hug}.

\begin{table}[h]
\begin{center}
\begin{tabular}{ccc}
\(\delta\) & \(\|x_1-x(\delta)\|\) & \(\|x_2-x(2\delta)\|\) \\
\hline
1/16 &$4.23\times 10^{-4}$  & $4.87\times 10^{-5}$\\
1/32 & $1.15\times 10^{-4}$ & $6.56\times 10^{-6}$\\
1/64 & $3.00\times 10^{-4} $& $8.50\times 10^{-7}$\\
1/128 &$7.62\times 10^{-6}$ & $1.08\times 10^{-7}$\\
1/256&$1.93\times 10^{-6}$& $1.36\times 10^{-8}$
\end{tabular}
\end{center}
\caption{Errors in the \(x\) variable vs.\ \(\delta\) after one or two timesteps.}
\label{table:cancellation}
\end{table}

\begin{theorem}\label{th:convergence} \emph{(Second-order convergence.)} Fix \(x_0\), \(v_0\) and a time interval \([0,T]\) and consider the solution \((x(t),v(t))\) of
\eqref{eq:odebis1}--\eqref{eq:odebis2} with \(x(0) = x_0\), \(v(0)=v_0\) and the iterates \((x_k,v_k)\) generated by \eqref{eq:alg1}--\eqref{eq:alg3}. Then, as \(\delta\rightarrow 0\),
\[
\max_{k\delta\leq T} \|x_{k}-x(k\delta)\| = \mathcal{O}(\delta^2),
\]
and
\[
\max_{k\delta\leq T} \|v_{k}-\big(\vpar(k\delta)+(-1)^k\vperp(k\delta)\big)\|  = \mathcal{O}(\delta^2).
\]
\end{theorem}

The constants implicit in the \(\mathcal{O}\) notation in the theorem vary with
\(f\), \(x_0\), \(v_0\) and \(T\) only.

Since the ODE preserves the value of \(f\) (Theorem~\ref{th:odebis}, Property 4.), the \(\mathcal{O}(\delta^2)\) deviations in the value of \(f\) for the \(x_k\), proved in Property 6.\ of Proposition~\ref{prop:alg}, match the second order convergence stated in Theorem~\ref{th:convergence}.

\subsection{Dynamics} \label{subsec:worked-out-example}
The appearence of the system \eqref{eq:odebis1}--\eqref{eq:odebis2} is certainly  unfriendly and it is not easy to guess the behaviour of its solutions. To get some insight, we illustrate its dynamics  in the case of two examples.

Assume first that initially \(v(0)\) is \emph{tangent}, i.e.\ \(\vperp(0) = 0\). Then  \eqref{eq:ode3} shows that, for all \(t\), \(\vperp(t)=0\) and therefore \(v(t)=\vpar(t)\). According to \eqref{eq:ode2}, the system reduces to:
\begin{eqnarray*}
\frac{d}{dt} x &=& \vpar,\\
\frac{d}{dt} \vpar &=&-\Npperp(x)[\vpar]\vpar.
\end{eqnarray*}
These are the equations of motion for a particle that is constrained to remain on the manifold \({\mathcal M}(f(x(0)))\) when there are no external forces (other than the normal force \(-\Npperp(x)[\vpar]\vpar\) exerted by the constraint)  \cite[Chapter 7]{Leimkuhler_Reich_2005}. The preservation of \(\|\vpar(t)\|^2\) is just the preservation of kinetic energy that takes place due to the absence of working forces.

As a second example we look at the case where \(n=2\), \(m=1\) and \(f(x) = -ax_{(1)}^2-bx_{(2)}^2\), \(a,b>0\), as in Example~\ref{ex:gaussian}. As pointed out there this would arise when sampling from a bivariate, anisotropic Gaussian target distribution. (Recall that for isotropic Gaussians all proposals are accepted.) There are four scalar equations in \eqref{eq:odebis1}--\eqref{eq:odebis2} and, using the first integrals \(f(x)\) and \(\|v\|^2\), the system is reducible to a two-dimensional one.
For simplicity we restrict the attention to solutions where \(f(x(0)) = -1\) (the general case may be retrieved from this by rescaling \(a\) and \(b\)). We may then parameterize the level set (ellipse) \({\mathcal M}(f(x(0)))\) in terms of an angular variable \(\phi\) as follows:
\[
x_{(1)}= a^{-1/2}\cos \phi,\qquad x_{(2)}= b^{-1/2}\sin\phi.
\]
The unit vector tangent to the ellipse is
\[
{\bf t}(\phi) = \mu^{-1/2} [-b^{1/2} \sin \phi, a^{1/2}\cos\phi]^{\top},
\]
where
\begin{equation}\label{eq:mu}
\mu= a\cos^2\phi+b\sin^2\phi,
\end{equation}
and the corresponding unit normal vector is
\[
{\bf n}(\phi) = \mu^{-1/2} [ a^{1/2}\cos\phi, -b^{1/2} \sin \phi]^{\top}.
\]

We may parameterize \(\vpar\) and \(\vperp\) in terms of scalars \(p\) and \(n\):
\[
\vpar = p{\bf t},\qquad \vperp = n {\bf n}.
\]
In addition \(p\) and \(n\) are linked by the first integral \(p^2+n^2 = p(0)^2+n(0)^2\) established  in Part 3 of  Theorem~\ref{th:odebis}. Substitution of these expressions into \eqref{eq:odebis1}, \eqref{eq:odebis3}, leads, using the expressions for \(\Nppar\) and \(\Npperp\) given in Example~\ref{ex:gaussian} and after considerable algebra, to the system
\begin{eqnarray}\label{eq:edophi}
\frac{d}{dt} \phi& = &(ab)^{1/2} \mu^{-1/2} p,\\
\label{eq:edop}
\frac{d}{dt} p& =& \frac{1}{2} (ab)^{1/2} (a-b) \mu^{-3/2} \sin(2\phi) (c^2-p^2),
\end{eqnarray}
where \(\mu\) is the function of \(\phi\) in \eqref{eq:mu} and \(c^2 = p(0)^2+n(0)^2\) is \(\|v(0)\|^2\), \(c \geq 0\). The system has to be considered only in the strip \(|p|\leq c\) because \(n(0)\) is real. The  lines \(p=\pm c\) that bound the strip are invariant.
\begin{figure}[t]
	\centering
	\includegraphics[width=0.6\textwidth]{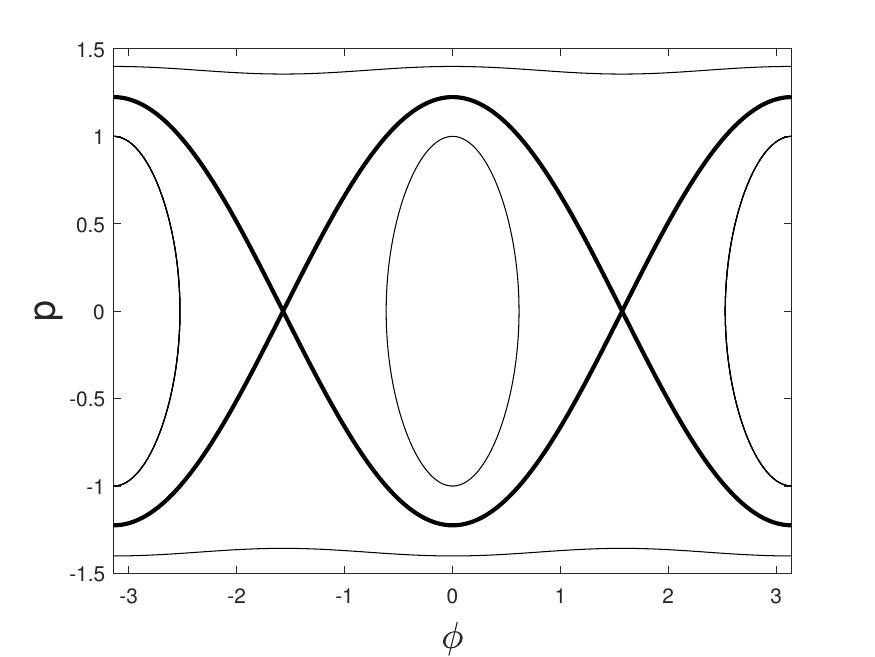}
	\caption{Phase plane of the system \eqref{eq:edophi}--\eqref{eq:edop}, when \(a=1\), \(b=4\), \(c=\sqrt{2}\). The portrait is \(2\pi\)-periodic in \(\phi\). The thicker lines correspond to the heteroclinic connections between saddle points. Outside the heteroclinic connections trajectories rotate; inside they librate.}
	\label{fig:phaseportrait}
\end{figure}

In the isotropic case \(a=b\), where the ellipses are circles, \(\mu\) does not vary with \(\phi\), \(p(t)\) remains constant and \(\phi(t)\) varies linearly with \(t\).  In the anisotropic case \(a\neq b\), it may be assumed that \(b>a\), as the other case is reduced to this by interchanging the roles of the components \(x_{(1)}\) and \(x_{(2)}\) of \(x\).
The system has equilibria at \(p=0\), \(\phi = j\pi/2\), \(j\) integer. One may easily prove that for \(j\) even the equilibria are centers and otherwise they are saddles.
A sketch of the phase portrait  may be seen in Figure~\ref{fig:phaseportrait}, where we note the similarity with the phase portrait of the standard pendulum equation. The saddle at  \(\phi = j\pi/2\), \(p=0\) (\(j\) odd) is connected to the neighbouring saddles at \(\phi = (j\pm2)\pi/2\), \(p=0\) via heteroclinic trajectories;
these are the image in the \((\phi,p)\)-plane of solutions \((\phi(t),p(t))\) with \[\lim_{t\rightarrow -\infty}(\phi(t),p(t))  = (j\pi/2,0),\qquad \lim_{t\rightarrow \infty}(\phi(t),p(t))  = ((j+2)\pi/2,0)  \]
or
\[\lim_{t\rightarrow -\infty}(\phi(t),p(t))  = (j\pi/2,0),\qquad \lim_{t\rightarrow \infty}(\phi(t),p(t))  = ((j-2)\pi/2,0).  \]

\begin{figure}[t]
	\centering
	\includegraphics[width=0.9\textwidth]{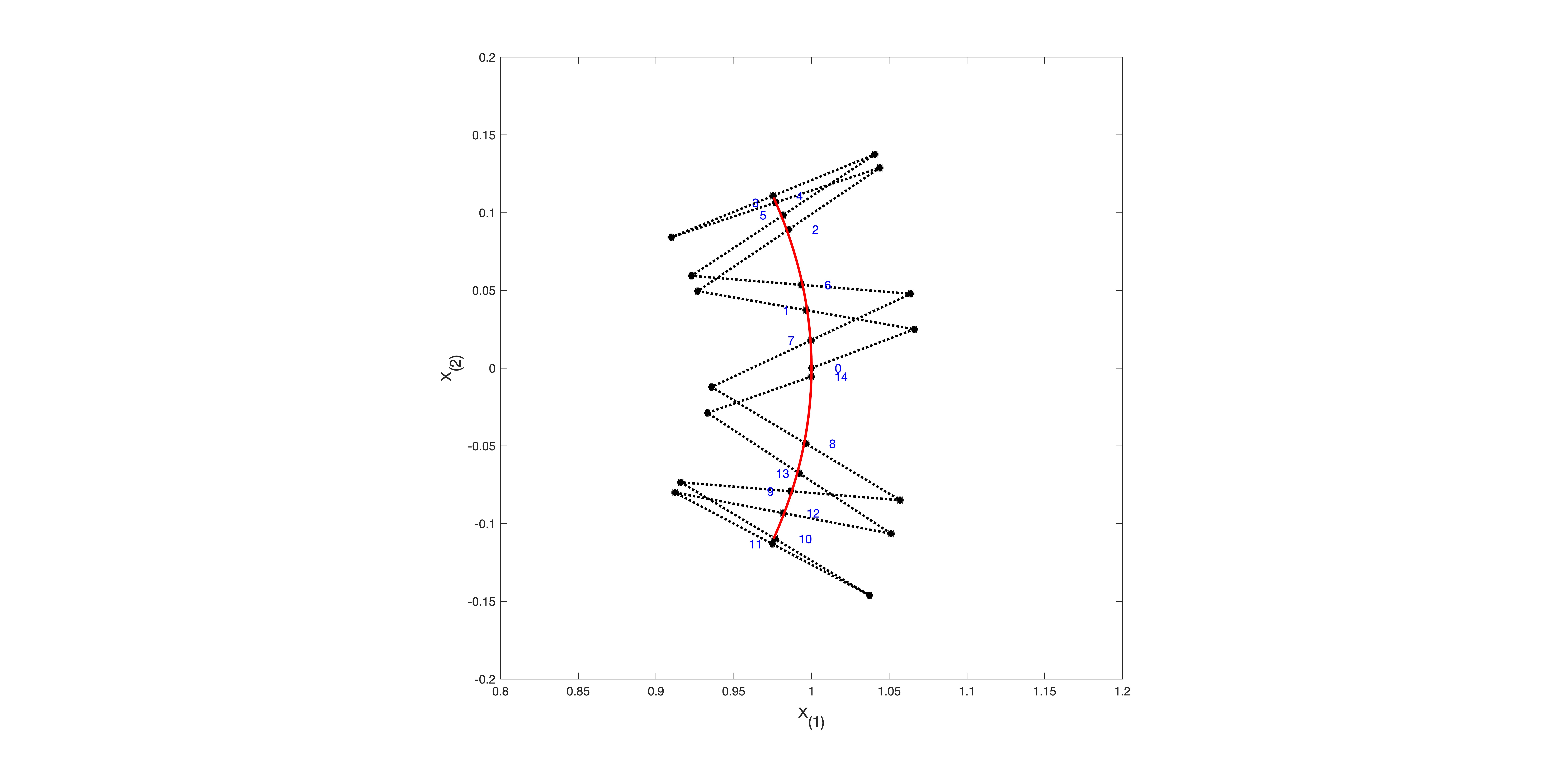}
	\caption{The continuous arc is the solution of the system \eqref{eq:odebis1}--\eqref{eq:odebis2} when \(f(x) = -x_{(1)}^2-4x_{(2)}^2\),
 \(x(0)=[1,0]\), and \(v(0)= [\sqrt{7/4}, 1/2]^{\top}\) depicted in the \(x\) plane for \(0\leq t \leq 1.4\). (This corresponds to system
\eqref{eq:edophi}--\eqref{eq:edop} with \(\phi(0) = 0\), \(p(0) = 1/2\), \(c=\sqrt{2}\).)
The solution \(x(t)\) first moves anticlockwise starting from \(x(0)\) but then folds back on itself, moves clockwise and folds back once more. At the final time, \(x(14)\) happens to be close to the initial location \(x(0)\). Also depicted is the approximation provided by Algorithm \eqref{eq:alg1}--\eqref{eq:alg3}, when \(\delta= 0.1\) and \(K=14\). The numerical solution mimics the behaviour of the ODE; the numbers \(0-14\) identify the stepnumber \(k\) of the iterates \(x_k\).}
	\label{fig:folding}
\end{figure}

 In the terminology of classical mechanics, trajectories at the top or the bottom of the figure (i.e.\ \lq outside\rq\ the heteroclinic connections) \emph{rotate}: in them \(p\) does not change sign as \(t\) varies  and the angle \(\phi\) increases monotonically (if \(p>0\)) or decreases monotonically (if \(p<0\)). The centers \(\phi = j\pi/2\), \(p=0\) (\(j\) even) are surrounded by trajectories that \emph{librate}: in them \(\phi\) varies periodically around \(j\pi/2\) between two values, say \(\phi_{\max},\phi_{\min}\), and \(p\) varies periodically around \(0\). In terms of the original system \eqref{eq:odebis1}--\eqref{eq:odebis2}, in rotating trajectories, as \(t\rightarrow\infty\), the point \(x(t)\) describes again and again the ellipse \({\mathcal M}(f(x(0)))\) either clockwise (if \(p(t)<0)\) or anticlockwise (if \(p(t)>0)\). For librating trajectories, the solution \(x(t)\)  describes the arc of the ellipse parameterized by \(\phi\in[\phi_{\max},\phi_{\min}]\); when \(\phi(t)\) has increased to \(\phi_{\max}\) it starts decreasing and the trajectory in the \(x\) plane folds back on itself, see Figure~\ref{fig:folding}. When \(\phi(t)\) has decreased to \(\phi_{\min}\) it starts increasing, causing the trajectory in the \(x\) plane to fold back on itself again. The process is repeated periodically.

 The figure also depicts the time-discrete solution  provided by  Hug, which is seen to mimic the behaviour of the ODE solution (this has to be the case, if \(\delta\) is sufficiently small, due to the convergence of the algorithm). Needless to say, this folding back phenomenon is not welcome in the context of Algorithm~\ref{alg:hug}, where one had the hope that, for suitable choices of \(K\), the \(x_k\) would wholly explore \({\mathcal M}(f(x_0))\).

Note from Figure~\ref{fig:phaseportrait} that, generally speaking, the (desirable) rotating trajectories correspond to cases where \(|p(0)|\) is large or, equivalently, \(|n(0)| = \sqrt{c^2-p(0)^2}\) is small relative to \(c=\|v(0)\|\). This could have been expected because for \(n(0) =0\) the system just describes the inertial motion of a particle around  the whole ellipse. In the scenario of Algorithm~\ref{alg:hug}, \(v(0)\) is a random variable with distribution \(q(\cdot|x(0))\). For the standard choice where \(q\) is an \(x(0)\)-independent isotropic normal, \(|n(0)|\) will not, on average, be small relative to \(\|v(0)\|\) with the result that, if the numerical trajectories are long, i.e.\ \(K\) is large,
many may fold back on themselves. For Algorithm~\ref{alg:hug}, this suggests (i) using distributions \(q(v|x)\) where at each \(x\) vectors with larger tangential components are given more weight and (ii) not using very large values of \(K\). We however expect this phenomenon to vanish as dimension increases, as discussed in the next section.

\section{Numerical experiments and future work}\label{sec:Numerical-experiments-and}

As discussed in the introduction, this work was initially motivated
by recent developments in the context of MCMC methods where the dynamics
studied in the present paper underpins the design of updates of Markov
chains aiming to explore, approximately, the contours of a probability
density, or more generally a manifold. Analysis and observations in
earlier sections point to potentially unappealing features of this dynamics.
In the simple scenario of a two-dimensional ellipse we have provided
a precise analysis of the dynamic, in particular characterising initial
velocities leading to either rotation or libration. However the two
dimensional scenario is particular in that trajectories that fail
to describe the full ellipse turn back on themselves exactly, resulting
in a periodic behaviour. Natural questions are therefore what form
this phenomenon takes in higher dimensions and whether it is even
noticeable when used in the (random) proposal mechanism of an MCMC
algorithm. Indeed, it is for example known that for a random vector
$v\sim\mathcal{U}(\mathbb{S}^{n})$ (where $\mathbb{S}^{n}$ is the
unit hypersphere) and any fixed unit vector $u\in\mathbb{R}^{n}$
then $u^{\top}v$ tends to be statistically small in the sense that
for $h\in[0,1]$
\[
\mathbb{P}(|u^{\top}v|\geq h)=\frac{\int_{h}^{1}[1-w^{2}]^{(n-2)/2}{\rm d}w}{\int_{0}^{1}[1-w^{2}]^{(n-2)/2}{\rm d}w}\,,
\]
confirming concentration on smaller values as $n$ increases. This
suggests that the negative phenomenon outlined earlier may perhaps naturally
vanish in high dimensional scenarios. A full, and relevant, analytical
investigation of these issues is beyond the scope of the present
manuscript and instead we present here a simple numerical study providing
some insight into this phenomenon in the multivariate scenario.

\begin{figure}
\centering
\subfloat[$\|v_{\perp}(0)\|=0.2023$]{\includegraphics[width=0.45\textwidth]{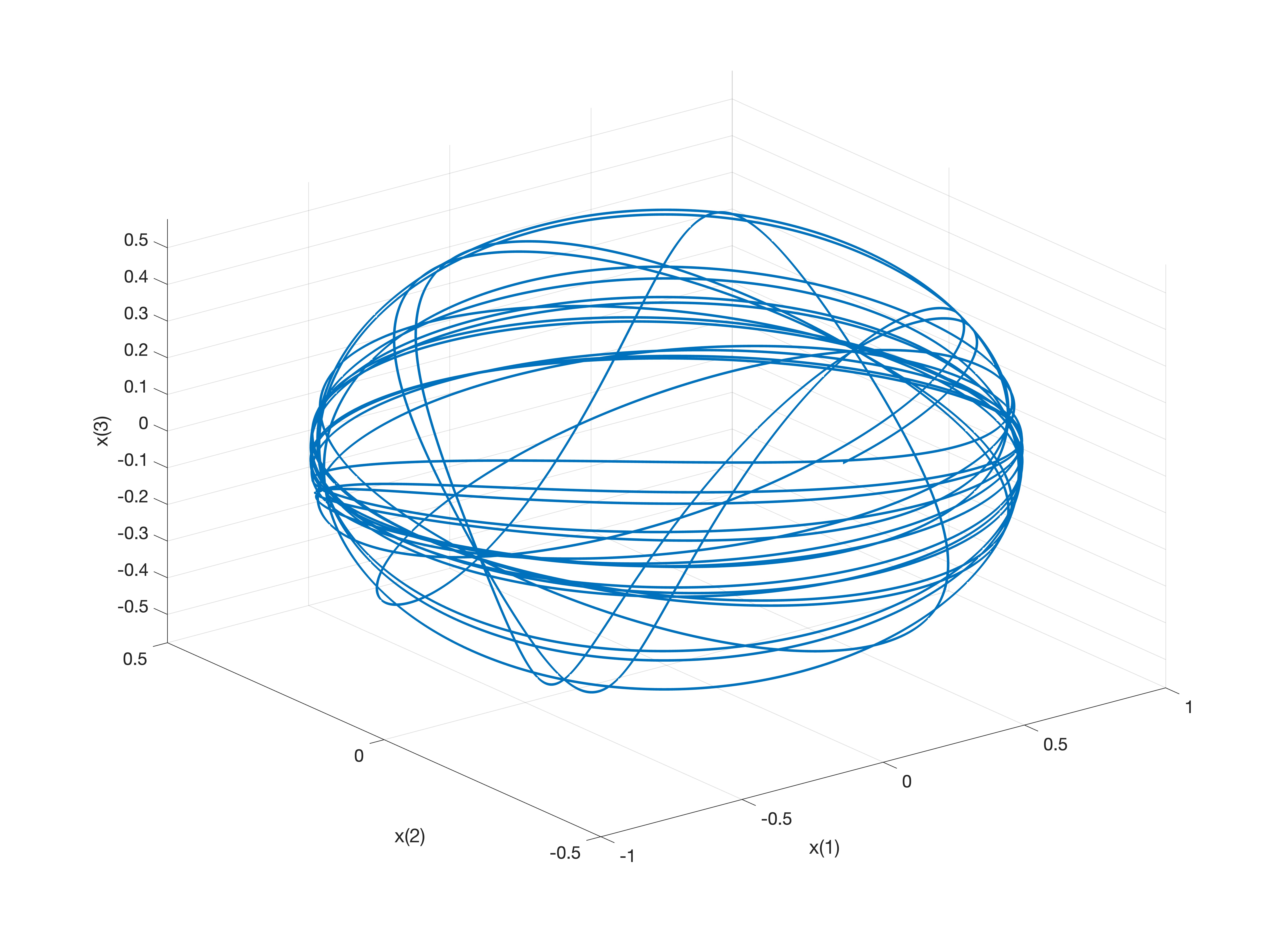}

}\subfloat[$\|v_{\perp}(0)\|=0.5673$]{\includegraphics[width=0.45\textwidth]{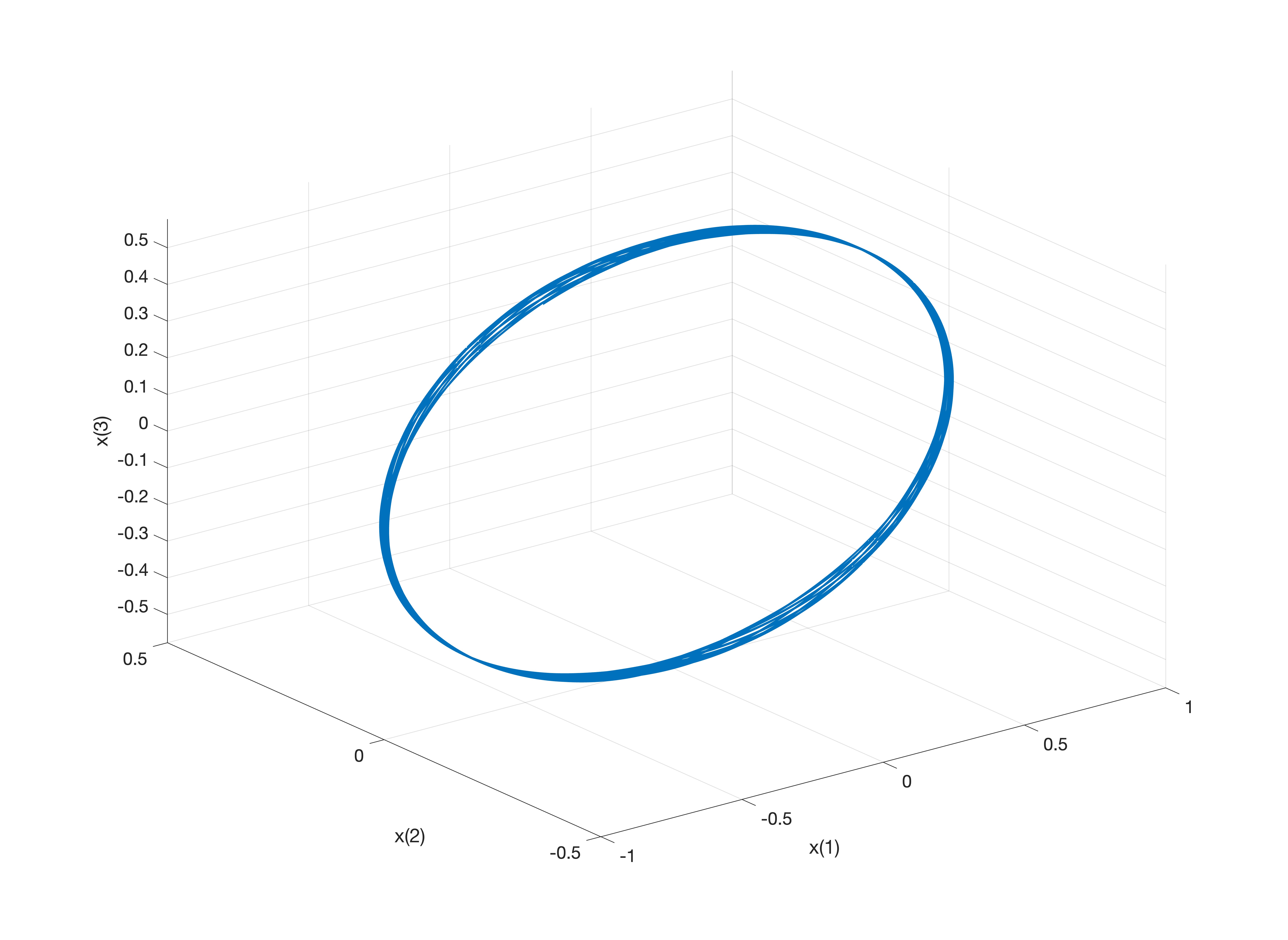}}

\subfloat[$\|v_{\perp}(0)\|=0.6647$]{\includegraphics[width=0.45\textwidth]{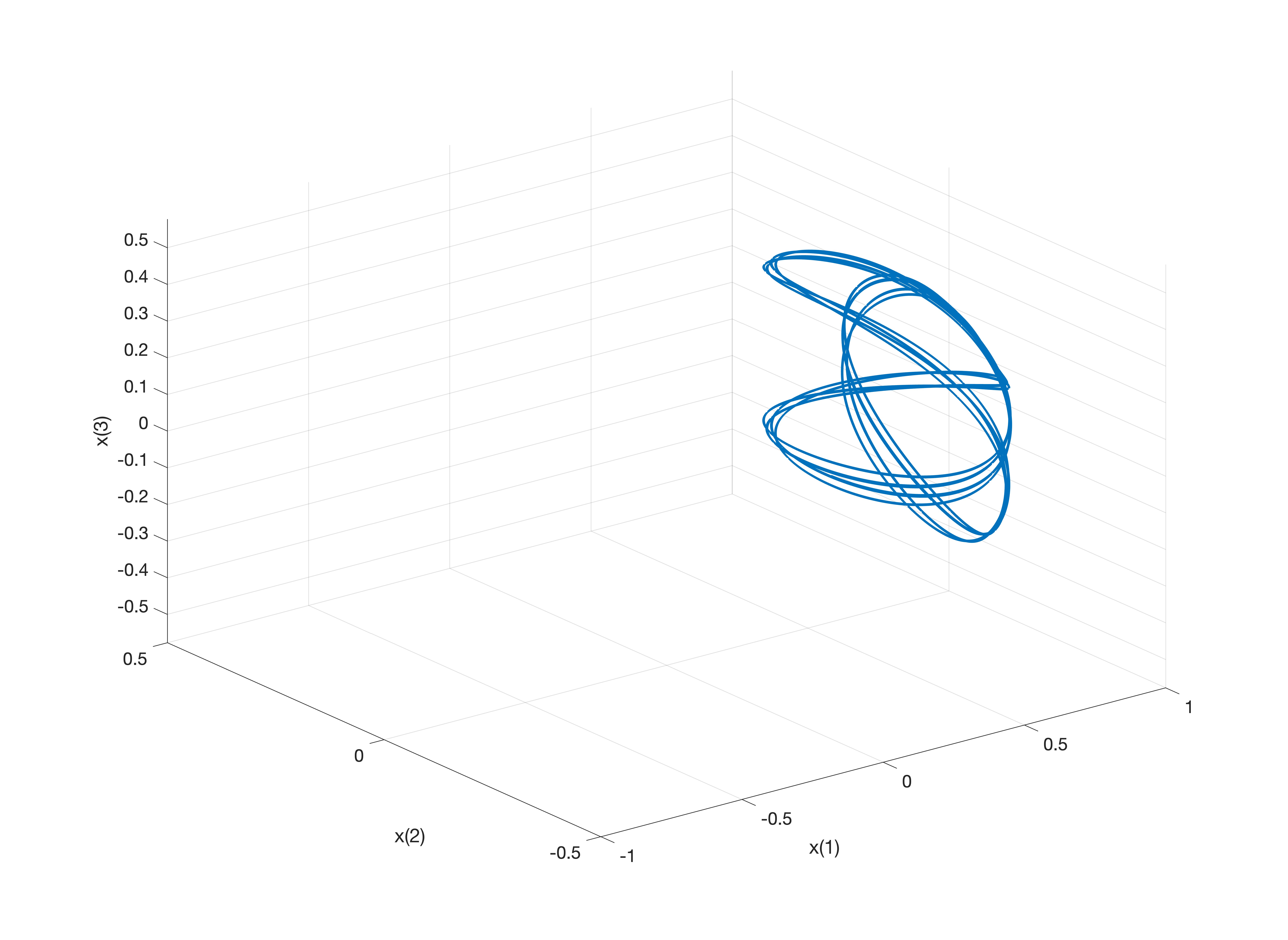}}\subfloat[$\|v_{\perp}(0)\|=0.7357$]{\includegraphics[width=0.45\textwidth]{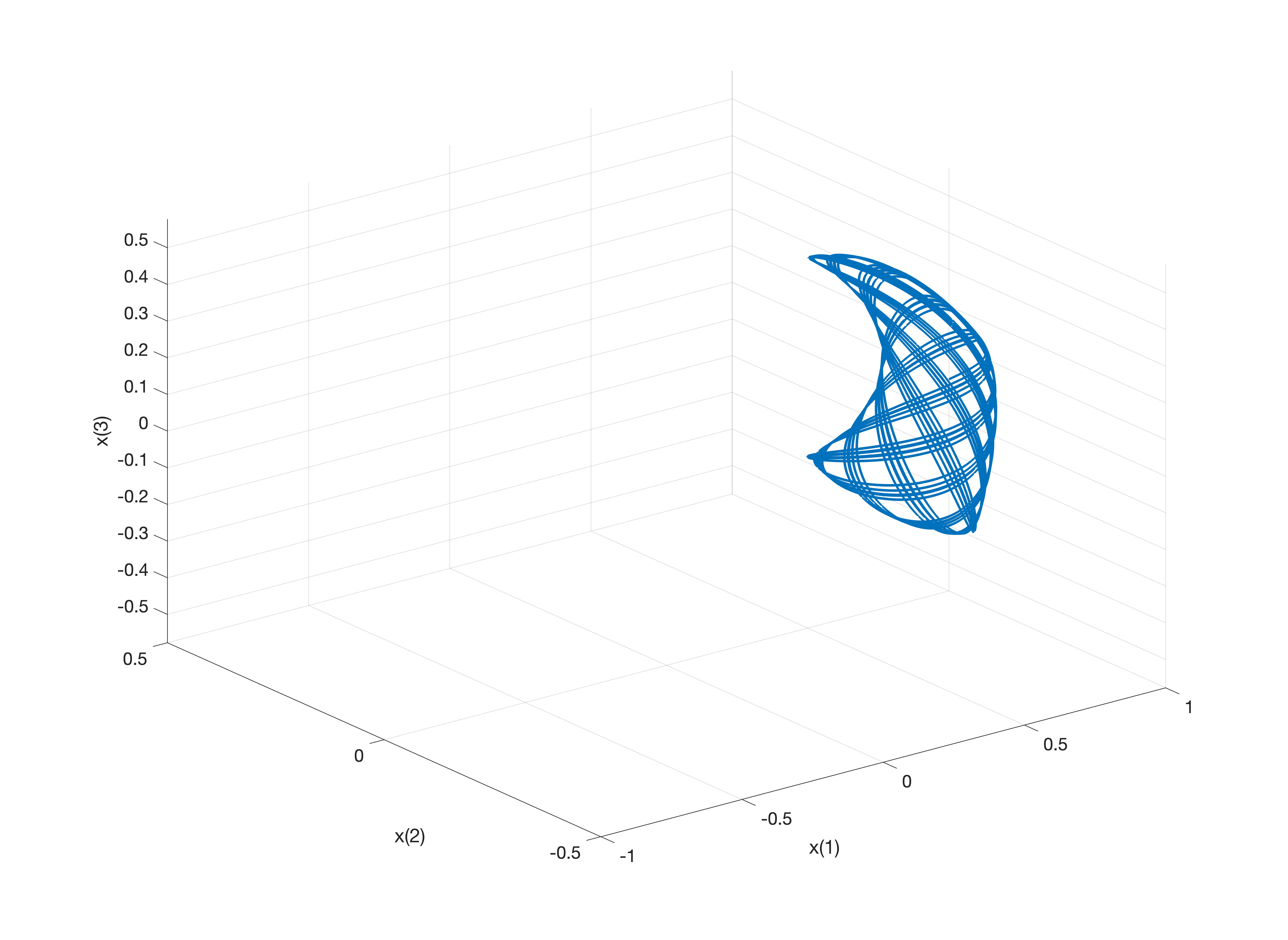}

}

\caption{$K=10000$ integration steps of the dynamics constrained to a $3D$
ellipsoid $(n=3)$, for four initial unit length velocities.}
\label{fig:3D-example-4-init-v}
\end{figure}

We first explore the behaviour of the Hug integrator in the scenario
where $n=3$ and the constraint of interest is an ellipsoid defined
by $f(x)=x^{T}Ax$ with $A:={\rm {\rm diag}}(1,4,3)$. In Fig.~\ref{fig:3D-example-4-init-v}
we display $K=10000$ steps of the integrator for $\delta=0.01$
and initial position $x(0)=[1,0,0]^{T}$ (therefore defining the level
set of $f$) for four initial velocities $v(0)$ such that $\|v(0)\|=1$.
These are representative of what we have observed for numerous draws
$v(0)\sim\mathcal{U}(\mathbb{S}_{n})$; as pointed out by one of the reviewers, $\|v_\perp(0)\|\sim \mathcal{U}(0, 1)$ in the specific scenario where $n=3$. As expected, we observe that a large orthogonal component
$\|v_{\perp}(0)\|$ results in trajectories that fail to explore a
large area of the ellipsoid, due to a near oscillatory phenomenon.
In fact we have found the size of the region visited to shrink for
even larger values $\|v_{\perp}(0)\|$, but do not report the results
here. For smaller values of $K$, more likely to be relevant in the
context of Monte Carlo algorithms, we have similarly observed the
influence of $\|v_{\perp}(0)\|$ on the size of the region explored
by the dynamics.

\begin{figure}
\centering
\subfloat[$K=1000$]{\includegraphics[totalheight=0.27\textheight]{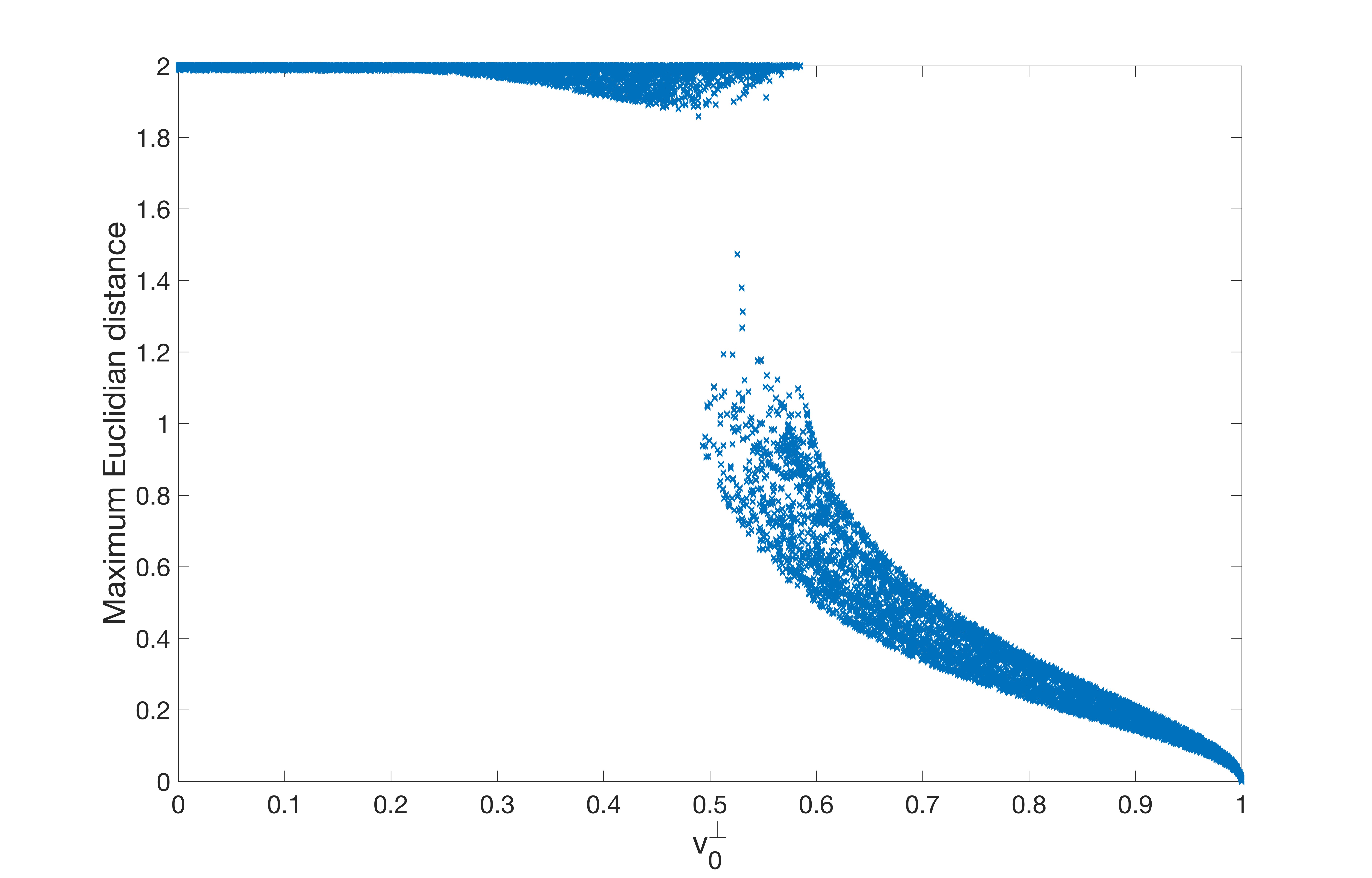}

}

\subfloat[$K=100$]{\includegraphics[totalheight=0.27\textheight]{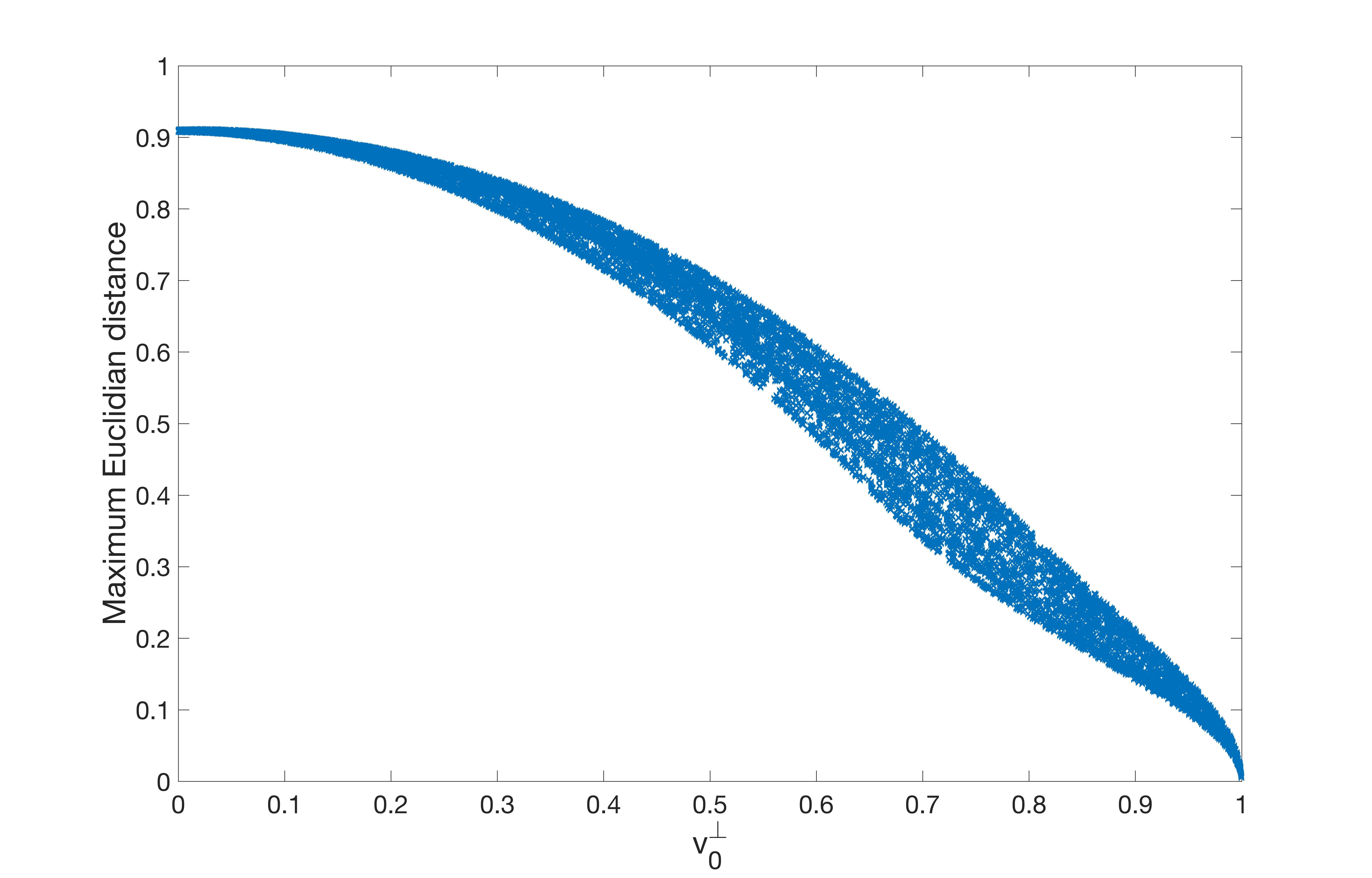}}

\subfloat[$K=10$]{\includegraphics[totalheight=0.27\textheight]{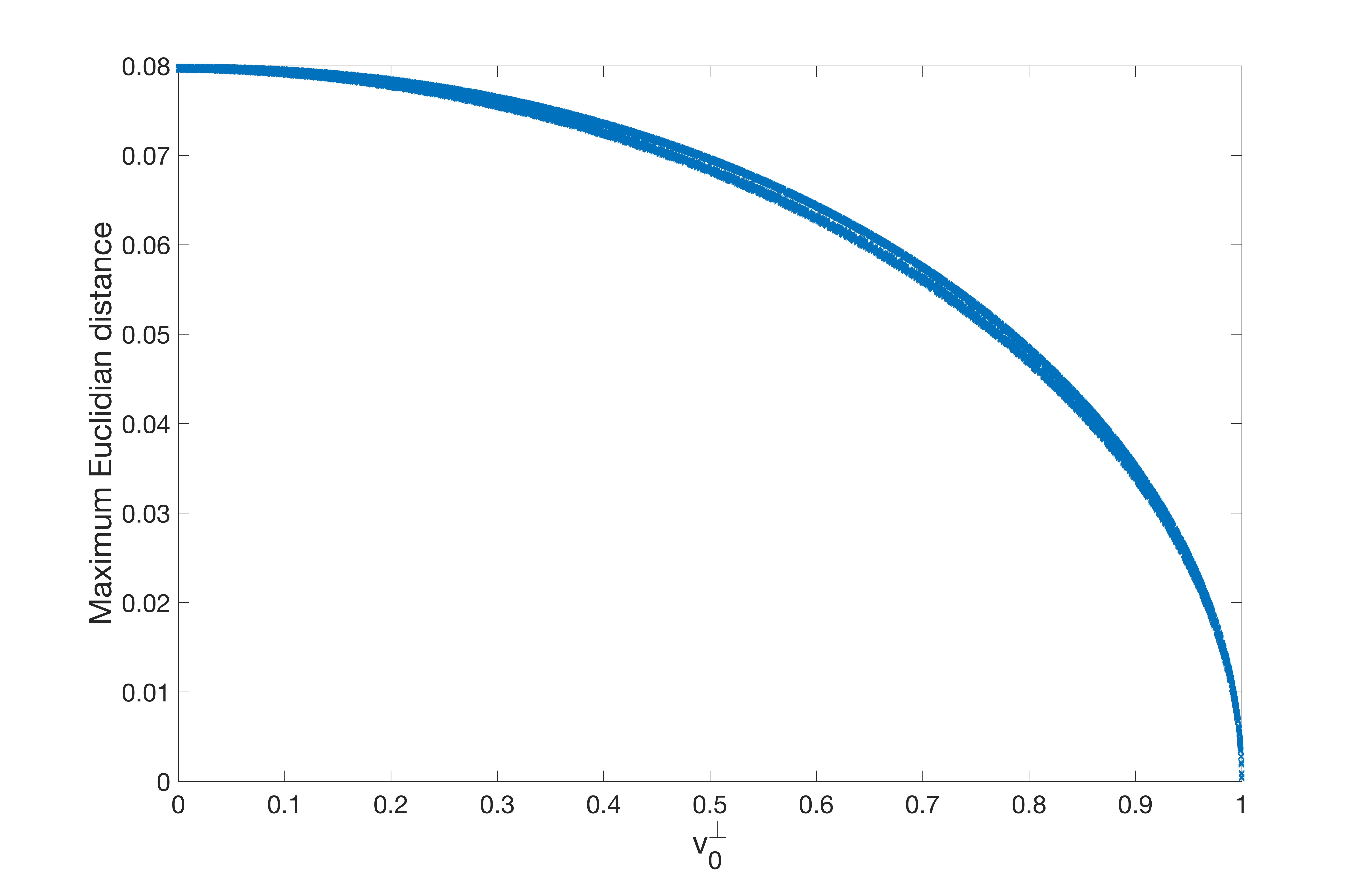}}\caption{Scatter plots $i\protect\mapsto\big(\|v_{\perp}^{(i)}(0)\|,d_{\rm max}\big(v_{\perp}(0)\big):=\max_{0\protect\leq k\protect\leq K}\|x^{(i)}(k)-x(0)\|\big)$
($n=3$)}
\label{fig:scatter-plots-3D}

\end{figure}

\begin{figure}
\centering
\subfloat[$K=1000$]{\includegraphics[totalheight=0.27\textheight]{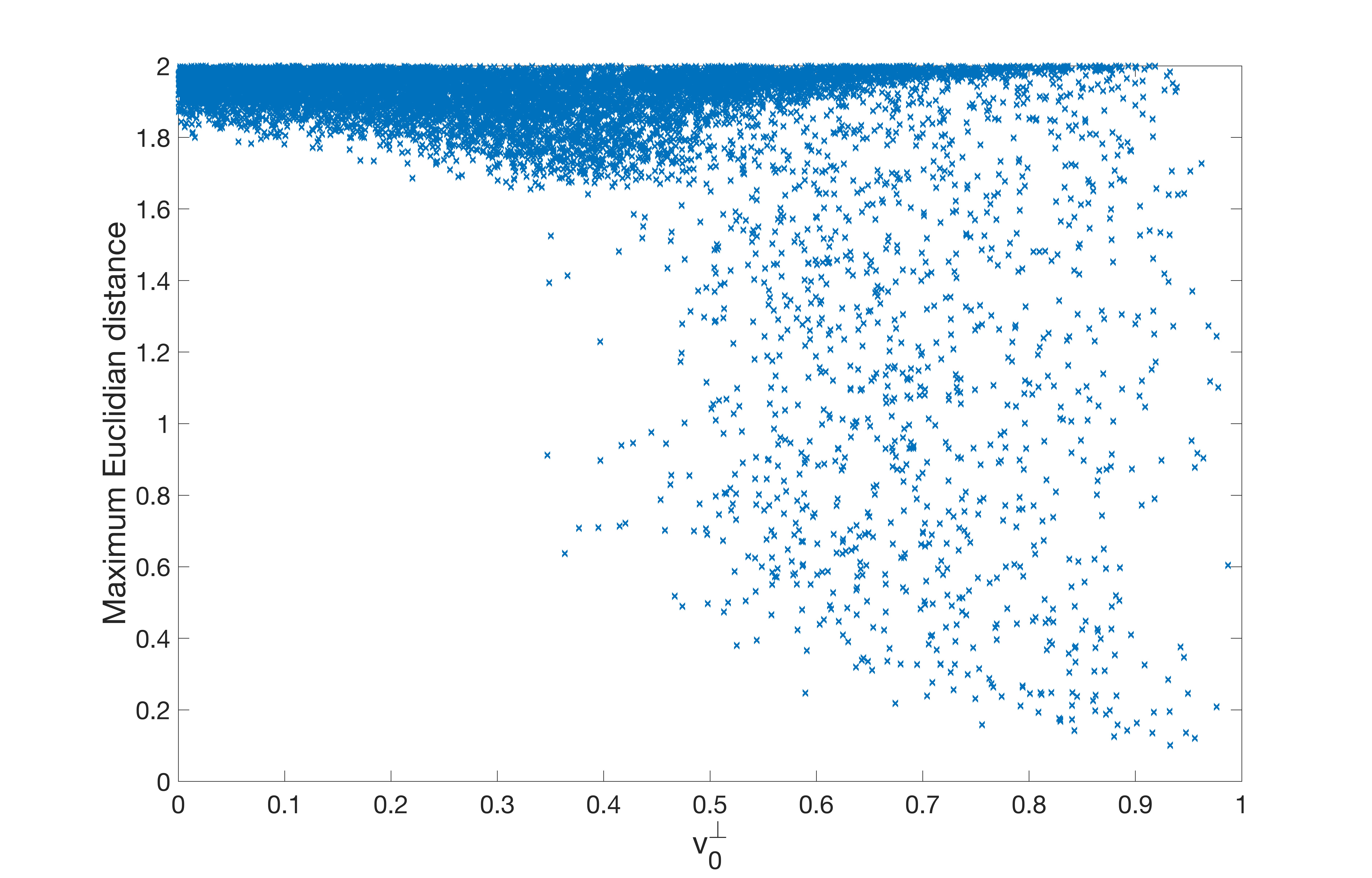}

}

\subfloat[$K=100$]{\includegraphics[totalheight=0.27\textheight]{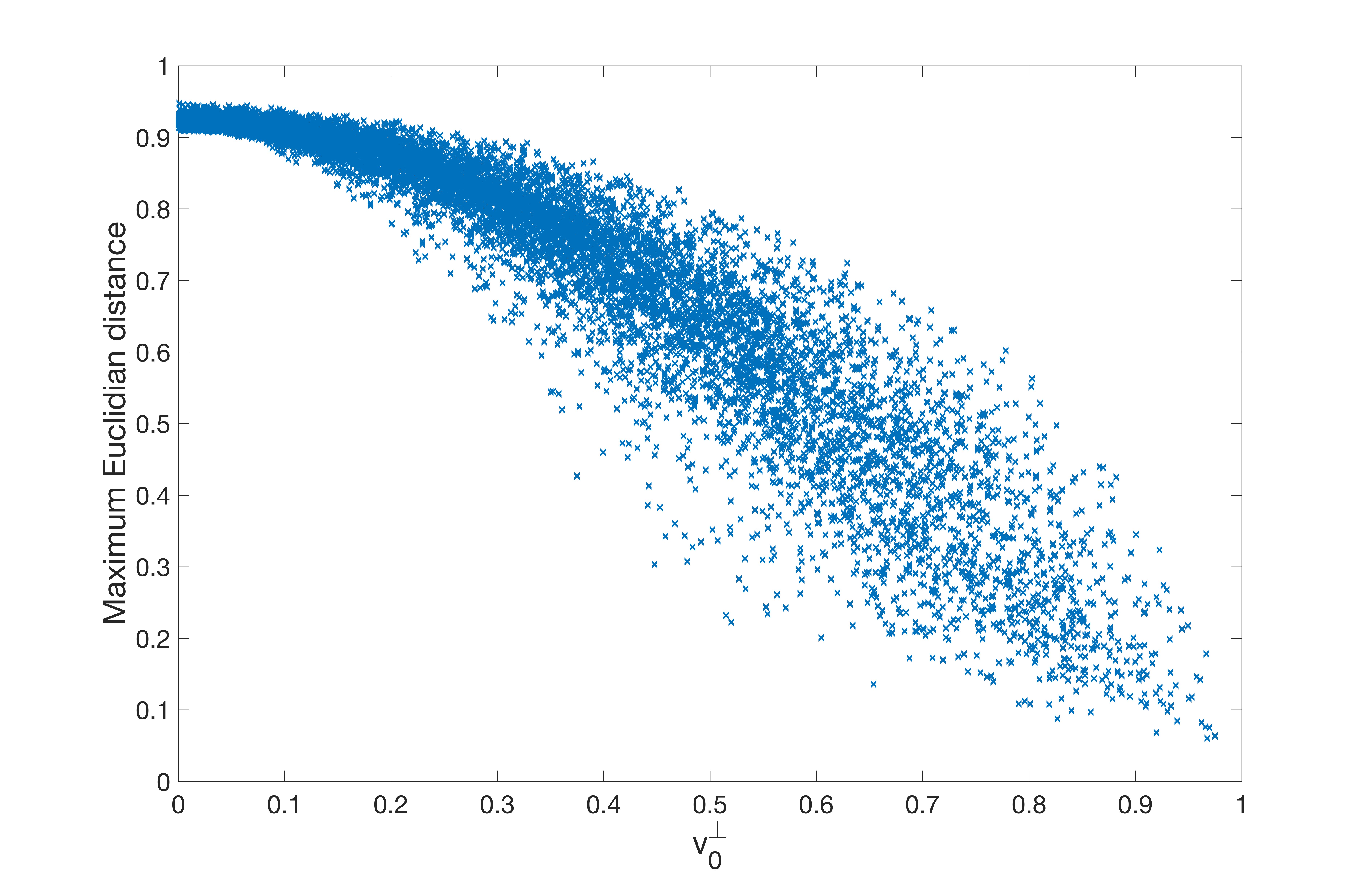}}

\subfloat[$K=10$]{\includegraphics[totalheight=0.27\textheight]{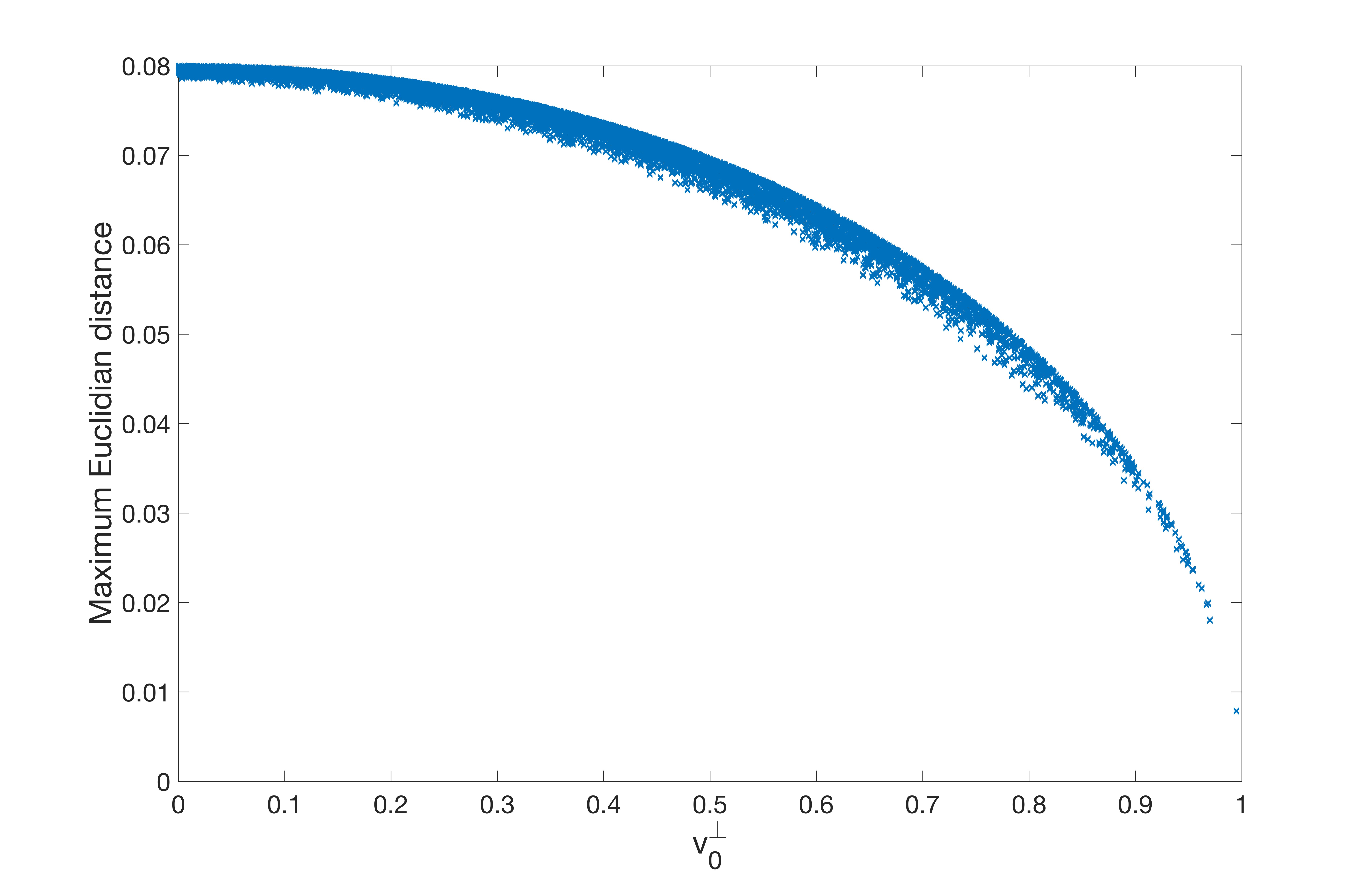}}

\caption{Scatter plots $i\protect\mapsto\big(\|v_{\perp}^{(i)}(0)\|,d_{\rm max}\big(v_{\perp}(0)\big):=\max_{0\protect\leq k\protect\leq K}\|x^{(i)}(k)-x(0)\|\big)$
($n=6$)}
\label{fig:6D-scatter}
\end{figure}

\begin{figure}
\centering
\subfloat[$K=1000$]{\includegraphics[totalheight=0.27\textheight]{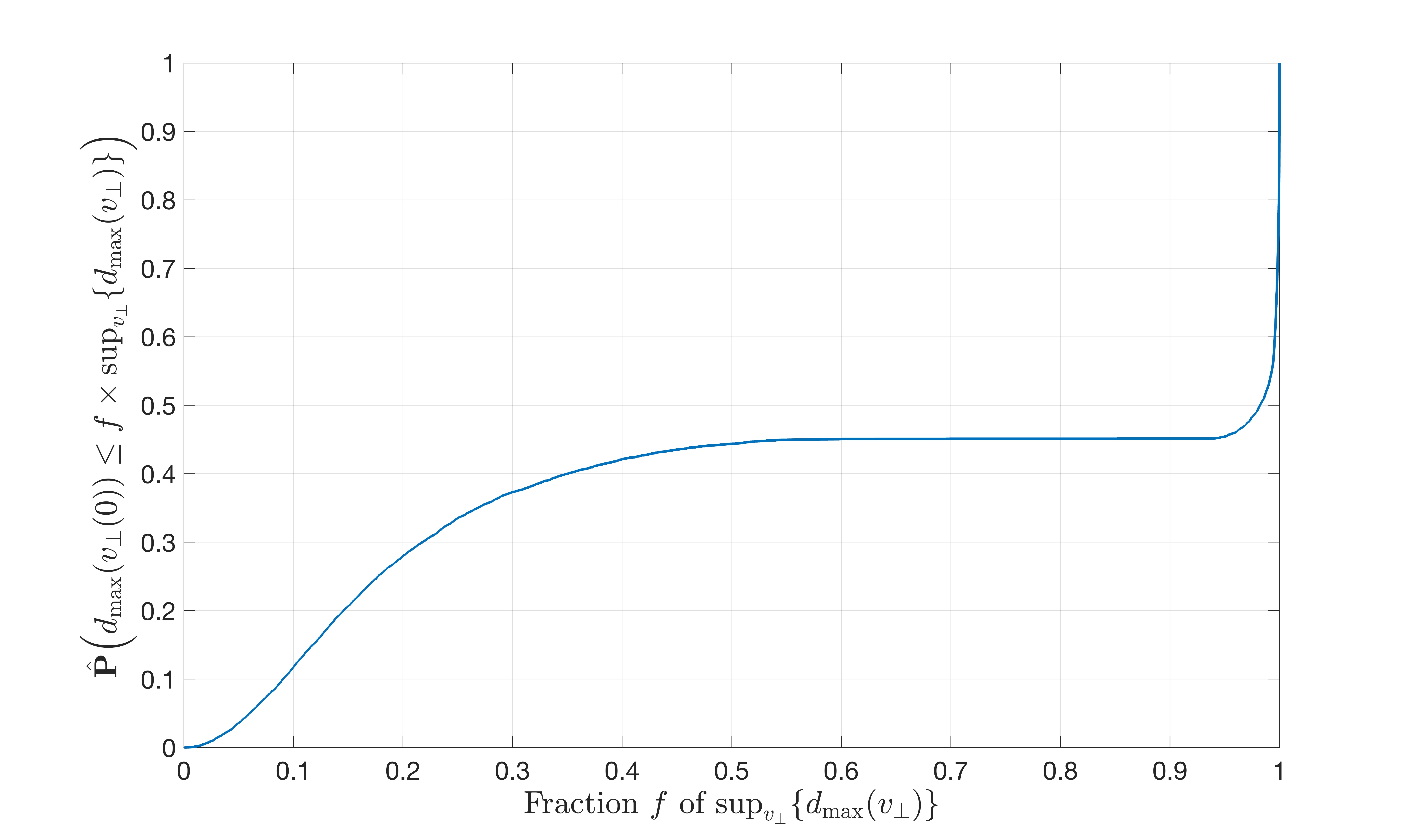}

}

\subfloat[$K=100$]{\includegraphics[totalheight=0.27\textheight]{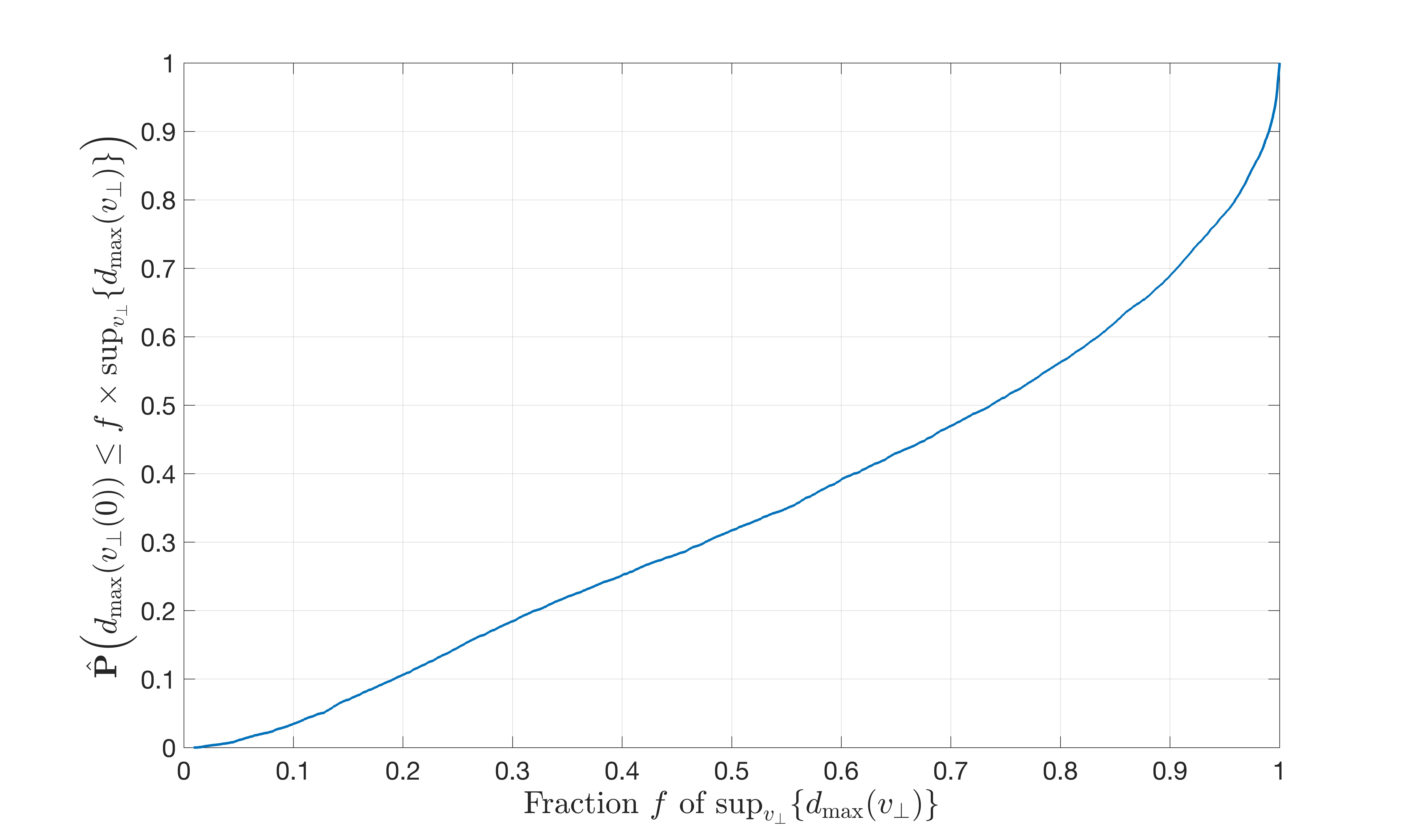}}

\subfloat[$K=10$]{\includegraphics[totalheight=0.27\textheight]{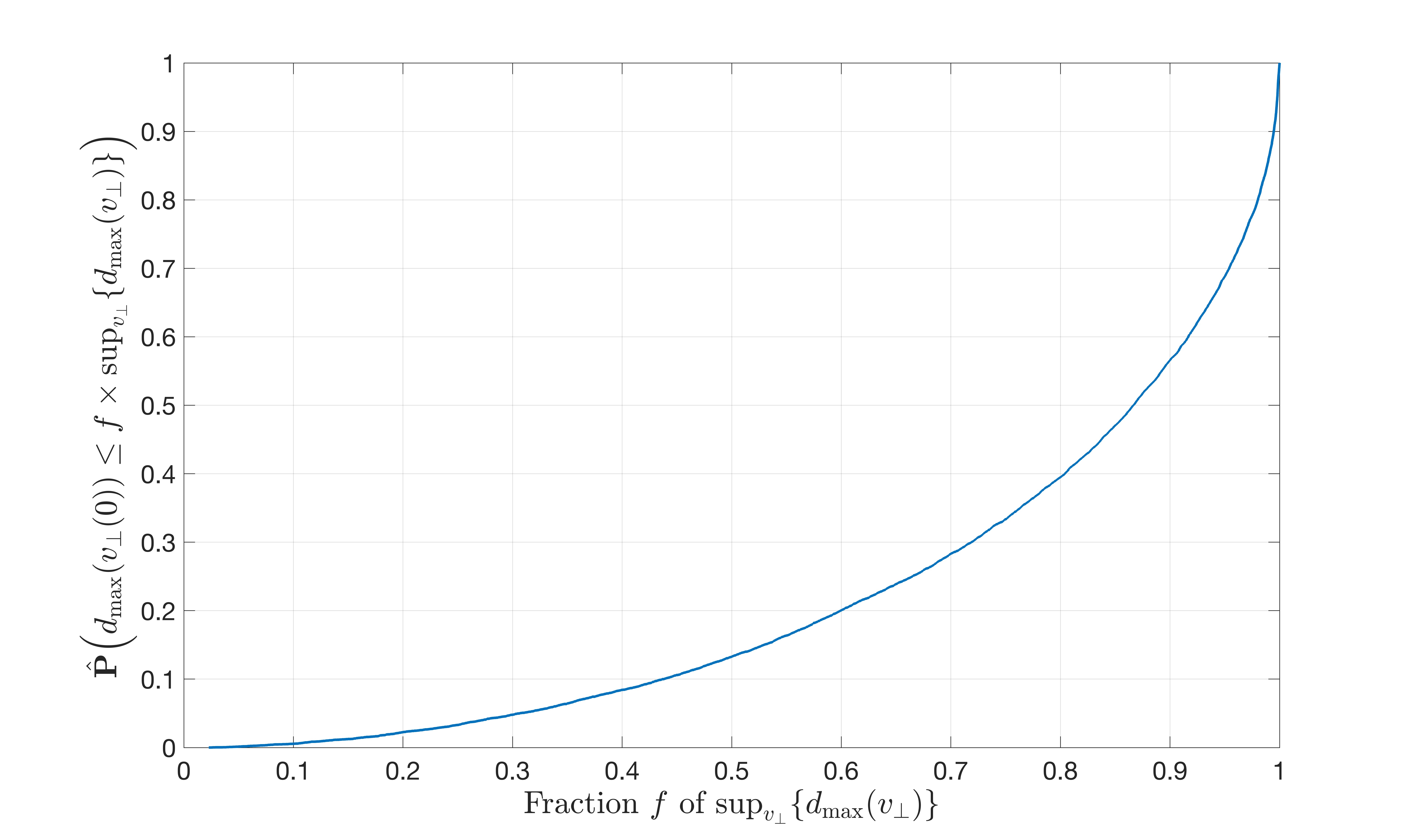}}

\caption{Estimates of $f\mapsto \mathbb{P}\big(d_{\rm max}\big(v_{\perp}(0)\big) \leq f \cdot \sup_{v_{\perp}} d_{\rm max}\big(v_{\perp}\big)\big)$ for $n=3$.}
\label{fig:3D-cdf}

\end{figure}

\begin{figure}
\centering
\subfloat[$K=1000$]{\includegraphics[totalheight=0.27\textheight]{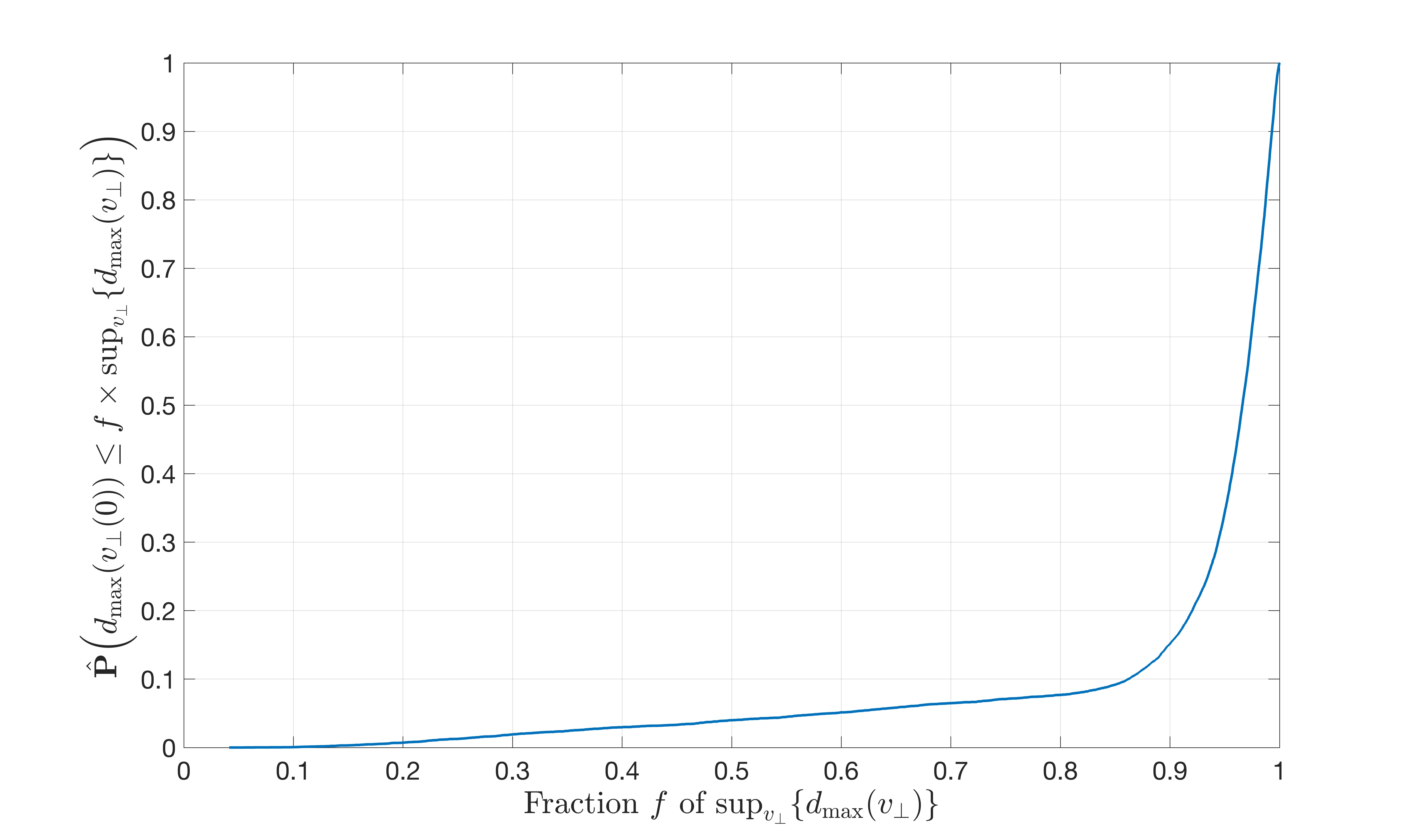}

}

\subfloat[$K=100$]{\includegraphics[totalheight=0.27\textheight]{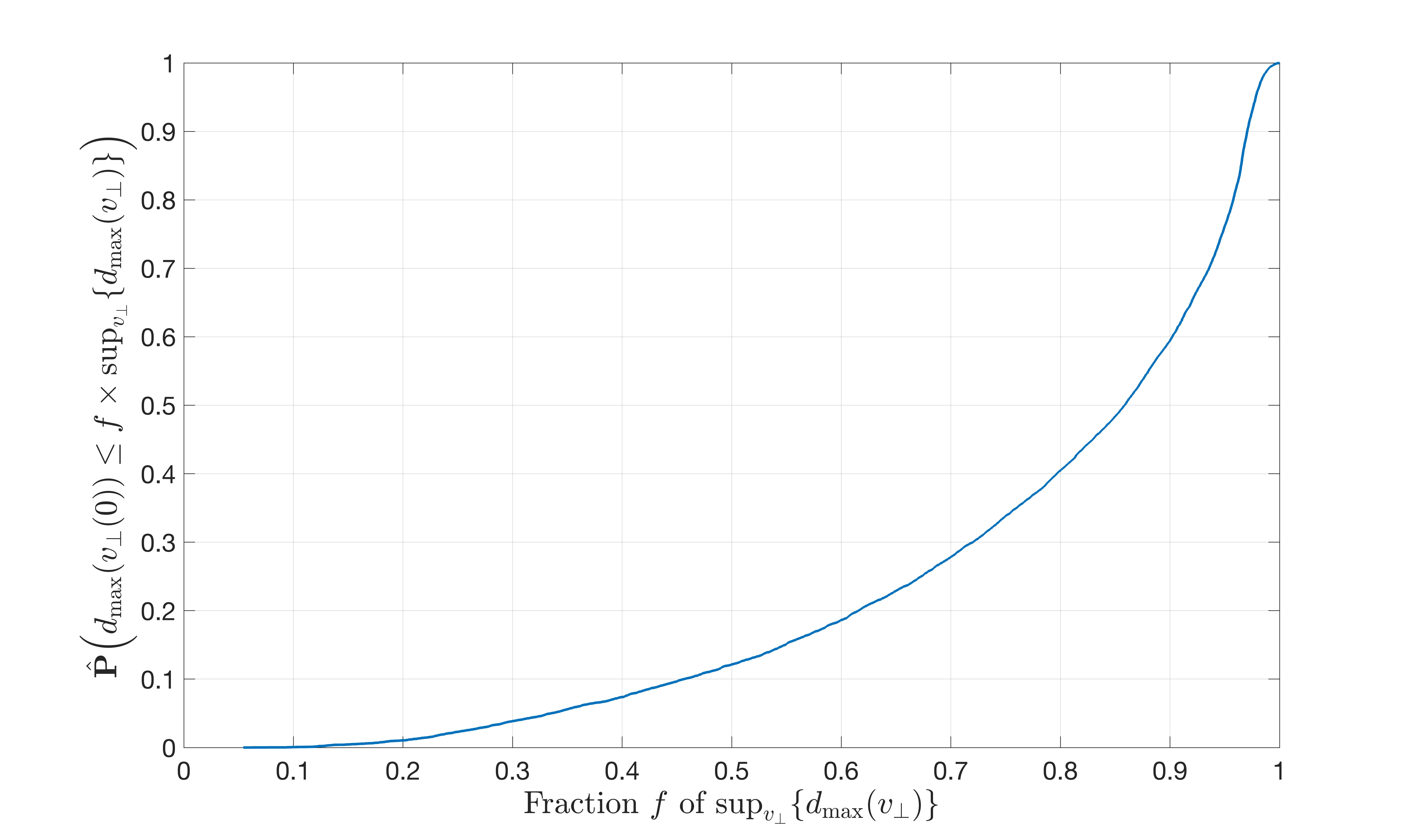}}

\subfloat[$K=10$]{\includegraphics[totalheight=0.27\textheight]{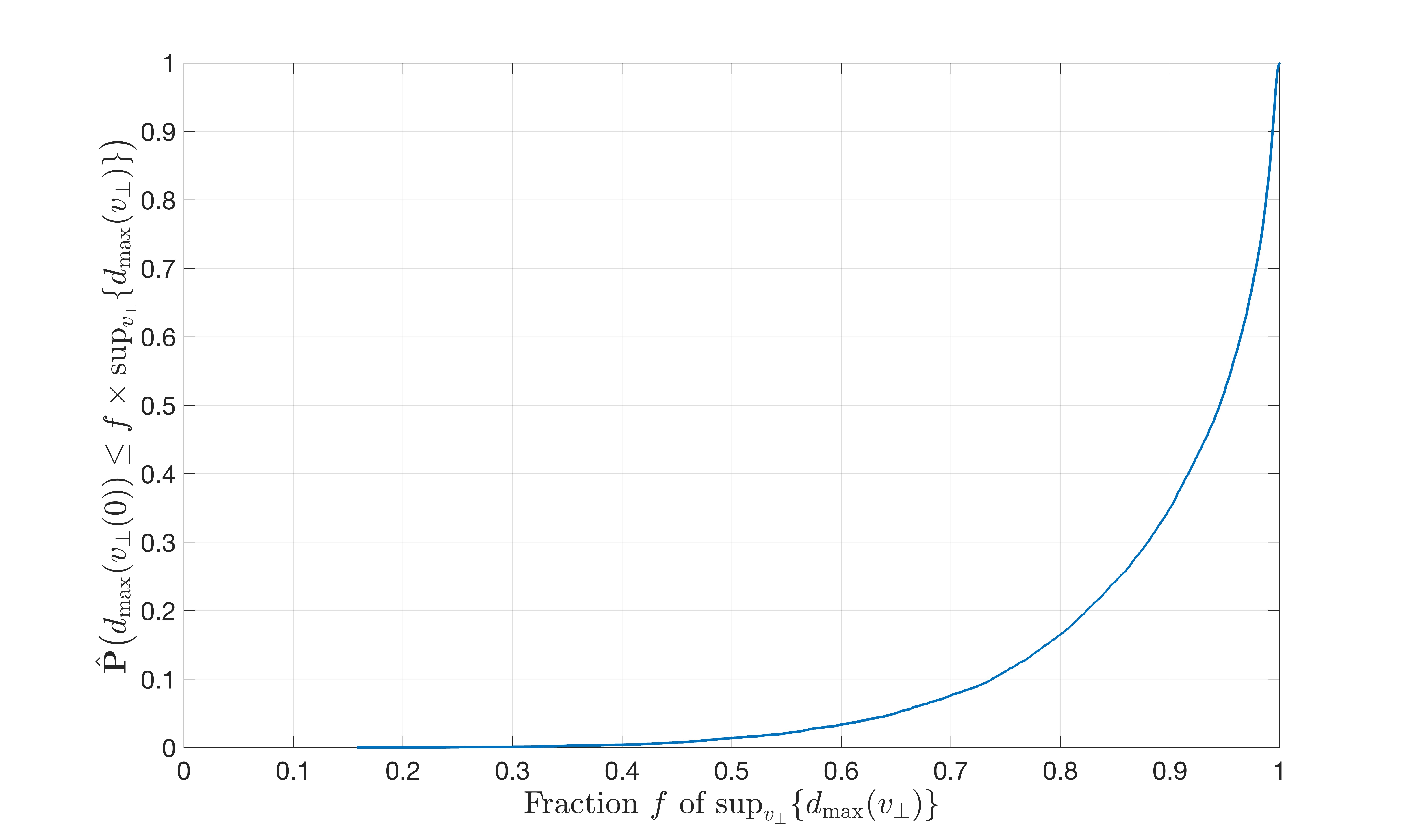}}

\caption{Estimates of $f\mapsto \mathbb{P}\big(d_{\rm max}\big(v_{\perp}(0)\big) \leq f \cdot \sup_{v_{\perp}} d_{\rm max}\big(v_{\perp}\big)\big)$ for  $n=6$.}
\label{fig:6D-cdf}
\end{figure}

In order to gain more insights into this phenomenon, we have run the
following experiments. In the setup above we now consider $10000$
realisations of the initial velocity $v^{(i)}(0)\overset{{\rm iid}}{\sim}\mathcal{U}(\mathbb{S}_{n})$
for $i=1,\ldots,10000$ and report the scatter plot $i\mapsto\big(\|v_{\perp}^{(i)}(0)\|,d_{\rm max}\big(v_{\perp}(0)\big):=\max_{0\leq k\leq K}\|x^{(i)}(k)-x(0)\|\big)$,
that is explore dependence of the largest Euclidean distance reached
by the dynamic on $\|v_{\perp}(0)\|$. The results are reported in
Fig.~\ref{fig:scatter-plots-3D} for $K=10,100,1000$ integration
steps. This confirms our earlier observations, including the concentration
of the whole dynamics around the initial position as $\|v_{\perp}(0)\|$
increases. We repeated this experiment for $n=6$, $A={\rm {\rm diag}}(1,4,3,5,1,10)$
and $x(0)=(1,0,0,0,0,0)^{T}$ and the results are reported in Fig.~\ref{fig:6D-scatter}
with similar conclusions. However, following a reviewer's suggestion, we present in Fig~ \ref{fig:3D-cdf} and Fig~ \ref{fig:6D-cdf} the Empirical Cumulative Distribution Functions (ECDF) of the maximum distance achieved for each realisations as a fraction of the maximum over all realisations, effectively approximating $f\mapsto \mathbb{P}\big(d_{\rm max}\big(v_{\perp}(0)\big) \leq f \cdot \sup_{v_{\perp}} d_{\rm max}\big(v_{\perp}\big)\big)$. As expected, for the reasons outlined at the beginning of the section, these ECDF functions are more skewed to higher fraction values for $n=6$ dimensions than for $n=3$, although we have not investigated the effect of the different extreme curvatures on the results.

Our work suggests that MCMC based on Hug trajectories may benefit
from the use of state-dependent and non-isotropic distributions for
the velocity, limiting the likelihood of using velocities with dominant
$v_{\perp}$ component; this however introduces a necessary tradeoff
since the dynamics now does not leave the velocity distribution invariant,
resulting in a higher rejection rate in the MH update \cite{C24}.
More work is needed to understand this tradeoff.

Note that, in cases where Hug is used to describe a manifold (as in Example~\ref{ex:splitting}), it is  possible to avoid the fold-back phenomenon by refreshing the velocity whenever its normal component has become large.

\section{Proofs}\label{sec:proofs}
This section contains the proofs of the results.

\subsection{Proof of Proposition~\ref{prop:alg}} \label{subsec:proof-prop-1}

We only have to establish Properties 5 and 6. With an obvious simplified notation, Taylor's theorem around \(x_{k+1/2}\) leads to
\begin{eqnarray*}
f_k & =& f_{k+1/2}-\frac{\delta}{2}J_{k+1/2}v_k +\frac{\delta^2}{4} H_k^+[v_k,v_k],\\
f_{k+1} & =& f_{k+1/2}+\frac{\delta}{2}J_{k+1/2}v_{k+1} +\frac{\delta^2}{4} H_{k+1}^-[v_{k+1},v_{k+1}],
\end{eqnarray*}
with
\begin{eqnarray*}
H_k^+ &= & \int_0^1 H\big(x_{k+1/2}+\lambda (x_k-x_{k+1/2})\big)(1-\lambda)d\lambda,\\
 H_{k+1}^-&=& \int_0^1 H\big(x_{k+1/2}+\lambda (x_{k+1}-x_{k+1/2})\big)(1-\lambda)d\lambda.
\end{eqnarray*}
Subtracting, taking into account that \(v_k+v_{k+1}\) is tangential by Property 3.\ and \eqref{eq:orthogonal},
\begin{equation}
\label{eq:fincrement}f_{k+1}-f_k = \frac{\delta^2}{4}  \Big( H_{k+1}^-[v_{k+1},v_{k+1}]-H_k^+[v_k,v_k]\Big)
\end{equation}
and, summing by parts,
\[
f(x_K)-f(x_0) =  \frac{\delta^2}{4}  \Big(  H_K^-[v_K,v_K]- H_0^+[v_0,v_0]+\sum_{k=1}^{K-1} \big(H_k^- -H_k^+\big)[v_k,v_k]\Big).
\]
By definition of the operator norm and Property 4., the norm of each of the terms \( H_K^-[v_K,v_K]\) and \(H_0^+[v_0,v_0]\) is bounded by
 \[ \beta \int_0^1 (1-\lambda) d\lambda\: \|v_0\|^2 = \frac{\beta}{2}\|v_0\|^2.
\]
On the other hand, \(H_k^- -H_k^+\) equals
\[
 \int_0^1 \Big(H\big(x_{k-1/2}+\lambda (x_{k}-x_{k-1/2})\big)-H\big(x_{k+1/2}+\lambda (x_k-x_{k+1/2})\big)\Big)(1-\lambda)d\lambda,
\]
and, therefore, by the assumed Lispchitz continuity, and Property 5.
\[
\|H_k^- -H_k^+\| \leq \int_0^1 \gamma (1-\lambda)^2 \|x_{k-1/2}-x_{k+1/2}\|\:d\lambda = \frac{\gamma}{3} \delta \|v_0\|.
\]
After putting everything together, we get the bound in 6.

With the hypothesis in  7., Property 4.\ shows that the right hand-side of \eqref{eq:fincrement} vanishes, which concludes the proof. When \(m=1\) the manifolds \({\mathcal M}(\eta)\) are concentric hyperspheres, and the property \(f(x_{k+1})-f(x_k)=0\) is apparent from the the geometry of the construction of \(x_{k+1}\) from \(x_k\).

\subsection{Proof of Lemma~\ref{lemma:nprime}}

1. From \eqref{eq:nprimeone}, the image of \(\Npperp(x)[w]\)  is contained in the image of \(J(x)^+\), which, according to \eqref{eq:J+} is contained in the image of \(J(x)^{\top}\), i.e.\ in the orthogonal of the kernel \({\mathcal T}(x)\) of \(J(x)\) (see \eqref{eq:orthogonal}), i.e.\ in \({\mathcal N}(x)\).

2. From \eqref{eq:nprimeone}, the kernel of \(\Npperp(x)[w]\)  is contained in the kernel of \(T(x)\), i.e.\ in \({\mathcal N}(x)\).

The other two claims are proved similarly.

\subsection{Proof of Theorem~\ref{th:odebis}}

For Item 1.\ we have to prove that the vector field in \eqref{eq:odebis1}--\eqref{eq:odebis2} is divergence free. We first present some facts related to the computation of divergences.
 If \(z\in\bbR^d\mapsto g(z)\in\bbR^d\) is a vector field, its divergence is of course the trace of the Jacobian
\begin{equation}\label{eq:div}
{\rm div}_z(g(z)) = \sum_j \langle e_j, g^\prime(z)e_j\rangle,
\end{equation}
here the \(e_j\) are the unit coordinate vectors and \(\langle\cdot,\cdot\rangle\) refers to the standard inner product. In the particular case of a linear vector field \(z\mapsto Az\), this reduces to the trace of the matrix \(A\):
\[
{\rm div}_z(Az) = \sum_j \langle e_j, Ae_j\rangle.
\]
Since the trace is invariant by similarity, if \(\Omega\) is an arbitrary invertible matrix,
\begin{equation}\label{eq:divsimilarity}
{\rm div}_z(Az) = \sum_j \langle e_j, Ae_j\rangle =
\sum_j \langle e_j, \Omega^{-1}A\Omega e_j\rangle.
\end{equation}
For a bilinear vector field \(z\mapsto B(z,z)\), we observe that the Jacobian is formed adding the Jacobian with respect to the first argument of \(B\) (with the second frozen) to the Jacobian with respect to the second argument (with the first frozen). Therefore the divergence may be written as
\begin{equation}\label{eq:divbilinear}
{\rm div}_z(B(z,z)) =
\sum_j \langle e_j, B(e_j,z) \rangle+\sum_j \langle e_j, B(z,e_j) \rangle.
\end{equation}

We are now ready to compute the divergence of the system \eqref{eq:odebis1}--\eqref{eq:odebis2}.
From \eqref{eq:div} and \eqref{eq:T},
\begin{equation}\label{eq:divx}
{\rm div}_x (T(x)v )= \sum_{j=1}^n \langle e_j,T^\prime(x)[e_j]v\rangle= -\sum_{j=1}^n \langle e_j,N^\prime(x)[e_j]v\rangle,
\end{equation}
and, from \eqref{eq:divbilinear}, with a shortened notation,
\begin{eqnarray}
\nonumber
&&{\rm div}_v \Big(\Big(\Nppar(x)\big[(T(x)-N(x))v\big]-\Npperp(x)\big[(T(x)-N(x))v\big]\Big)v\Big)  = \\
\nonumber
&&\qquad\qquad \sum_{j=1}^n \langle e_j, \Big(\Nppar\big[(T-N)e_j\big]-\Npperp\big[(T-N)e_j\big]\Big)v\rangle+\\
\label{eq:divv}
&& \qquad\qquad \sum_{j=1}^n \langle e_j, \Big(\Nppar\big[(T-N)v\big]-\Npperp\big[(T-N)v\big]\Big)e_j\rangle.
\end{eqnarray}
The last summation vanishes because it is the trace of the matrix \(\Nppar\big[(T-N)v\big]-\Npperp\big[(T-N)v\big]\), which is skew-symmetric due to \eqref{eq:nprimetwo}.\footnote{In actual fact the traces of \(\Nppar\big[(T-N)v\big]\) and \(\Npperp\big[(T-N)v\big]\) are both \(0\). According to Lemma~\ref{lemma:nprime}, these matrices are nilpotent with \(\Nppar\big[(T-N)v\big]^2=0\) and \(\Npperp\big[(T-N)v\big]^2 = 0\) which implies that all eigenvalues of \(\Nppar\big[(T-N)v\big]\) and \(\Npperp\big[(T-N)v\big]\) vanish.  } The properties of the images in Lemma~\ref{lemma:nprime} show that the first summation in the right hand-side of \eqref{eq:divv} may be rewritten as
\[
\sum_{j=1}^n \langle e_j, \Big((T-N)\Nppar\big[(T-N)e_j\big]+(T-N)\Npperp\big[(T-N)e_j\big]\Big)v\rangle
\]
or, resorting to \eqref{eq:Nprime},
\[
\sum_{j=1}^n \langle e_j, (T-N)N^\prime\big[(T-N)e_j\big]v\rangle.
\]
Since \((T-N)^{-1} = T-N\), the formula \eqref{eq:divsimilarity} implies that this equals
\[
\sum_{j=1}^n \langle e_j, N^\prime[e_j]v\rangle,
\]
which cancels \eqref{eq:divx}. In this way volumen conservation has been established.

 The time reversibility in Item 2.\ is clear:
changing \(t\) into \(-t\) and \(v\) into \(-v\) reverses the sign of the left and right hand-sides of \eqref{eq:odebis1} and does not change
the left and right hand-sides of \eqref{eq:odebis2} as the right hand-side is quadratic in \(v\).

For Item 3., from \eqref{eq:odebis3}
\[
\frac{1}{2} \frac{d}{dt}\langle v,v\rangle = \langle v, \Nppar(x)[\vpar-\vperp]\vperp\rangle\:-\:\langle v,\Npperp(x)[\vpar-\vperp]\vpar\rangle,
\]
or, taking into account the properties of the images of \(\Nppar(x)[\vpar-\vperp]\) and \(\Npperp(x)[\vpar-\vperp]\) in Lemma~\ref{lemma:nprime},
\[
\frac{1}{2}\frac{d}{dt} \langle v,v\rangle = \langle \vpar, \Nppar(x)[\vpar-\vperp]\vperp\rangle\:-\:\langle \vperp,\Npperp(x)[\vpar-\vperp]\vpar\rangle.
\]
The right hand-side vanishes because the matrices involved are transposed of one another (see \eqref{eq:nprimetwo}).

Finally, for Item 4., \eqref{eq:odebis1} leads to
\[
\frac{d}{dt} f(x) = J(x) \vpar;
\]
the right hand-side vanishes according to \eqref{eq:orthogonal}.
\subsection{Proof of  of Proposition~\ref{prop:components}}

By using \eqref{eq:T}, \eqref{eq:odebis3}, \eqref{eq:odebis1} and \eqref{eq:Nprime} successively, we may write
\begin{eqnarray*}
\frac{d}{dt}( T v )& =& \Big(\frac{d}{dt}T\Big)v+T\frac{d}{dt}v\\
&=& -\Big(\frac{d}{dt} N\Big) v+T\Nppar[\vpar-\vperp]\vperp-T\Npperp [\vpar-\vperp]\vpar\\
&=& -\big(\Nppar[\vpar]+\Npperp[\vpar]\big) (\vpar+\vperp)\\
&&\qquad\qquad+T\Nppar[\vpar-\vperp]\vperp-T\Npperp [\vpar-\vperp]\vpar.
\end{eqnarray*}
Simplification with the help of Lemma~\ref{lemma:nprime} leads to \eqref{eq:ode2}

 Equation \eqref{eq:ode3} is derived in a similar fashion.

\subsection{Proof of Theorem~\ref{th:consistency}}

For  $\sigma$:
\[\sigma_{k+1}  =
X_{k+1}-X_k -\delta T(X_k)V_k-
\delta \big(T(X_k+(\delta/2)V_k)-T(X_k)\big)V_k
\]
Now, note that
\[
T(X_k+(\delta/2)V_k)-T(X_k)= \mathcal{O}(\delta),
\]
and that, in view of \eqref{eq:odebis1},
\[
X_{k+1}-X_k -\delta T(X_k)V_k =x((k+1)\delta)-x(k\delta) -\delta \vpar(k\delta) = \mathcal{O}(\delta^2).
\]

We turn to $\tau$. We write:
\begin{eqnarray*}
\tau_{k+1}  & =&V_{k+1}- \big(I-2N(X_k)\big)V_k
 + 2\big(N(X_k+(\delta/2)V_k)-N(X_k)\big)V_k\\
&=&V_{k+1}-T(X_k)V_k+N(X_k)V_k
 + 2\big(N(X_k+(\delta/2)V_k)-N(X_k)\big)V_k,
\end{eqnarray*}
and Taylor expanding
\[
\tau_{k+1} = V_{k+1} -T(X_k)V_k+N(X_k)V_k
 +\delta N^\prime(X_k)[V_k]V_k+\mathcal{O}(\delta^2),
\]
or, from \eqref{eq:definition} (assuming \(k\) even) and \eqref{eq:Nprime}, simplifying slightly the notation:
\begin{eqnarray*}
\tau_{k+1}& = &\vpar((k+1)\delta)-\vperp((k+1)\delta)-\vpar(k\delta)+\vperp(k\delta)\\
&&+\delta \Nppar[\vpar(k\delta)+\vperp(k\delta)](\vpar(k\delta)+\vperp(k\delta))\\&& +\delta \Npperp[\vpar(k\delta)+\vperp(k\delta)]
(\vpar(k\delta)+\vperp(k\delta)) +\mathcal{O}(\delta^2).
\end{eqnarray*}
We next recall \eqref{eq:ode2} and \eqref{eq:ode3} and  Lemma~\ref{lemma:nprime} to get, simplifying further the notation:
\begin{eqnarray*}
\tau_{k+1}& = &-\delta\Nppar[\vperp]\vperp -\delta\Npperp[\vpar]\vpar-\delta \Nppar[\vpar]\vperp-\delta\Npperp[\vperp]\vpar\\&&
+\delta \Nppar[\vpar+\vperp]\vperp +\delta \Npperp[\vpar+\vperp]
\vpar +\mathcal{O}(\delta^2)\\&=& \mathcal{O}(\delta^2).
\end{eqnarray*}

For \(k\) odd the technique of proof is the same.

\subsection{Proof of Theorem~\ref{th:convergence}}

Standard textbook proofs (see e.g.\ \cite[Chapter 2.3]{HNW93}), going all the way back to Henrici \cite{H}, of the convergence of consistent one-step integrators assume that the
map that advances the numerical solution over a single time-step is  of the form \(Id+\delta \Delta\), where the so-called increment function \(\Delta\) is assumed to be Lipschitz in the neighbourhood of the ODE solution being approximated.
Due to the flipping of the normal component of the velocity, that assumption does not hold for \eqref{eq:alg1}--\eqref{eq:alg3}.
 There is a second difficulty. The  truncation errors \eqref{eq:sigma}--\eqref{eq:tau} are \(\mathcal{O}(\delta^2)\) (first-order consistency) and the usual consistency+stability argument would lead to \(\mathcal{O}(\delta)\) global error bounds (first order convergence), while Theorem~\ref{th:convergence} claims \(\mathcal{O}(\delta^2)\) error bounds (second-order convergence).

In order to circumvent these difficulties,  we consider the numerical integrator \(\widehat{\mathcal I}\) such that one timestep of length \(2\delta\) with \(\widehat{\mathcal I}\) is the result
of taking successively two steps of length \(\delta\) with  \eqref{eq:alg1}--\eqref{eq:alg3}. In symbols, the one step map \(\widehat{\Psi}\) associated with \(\mathcal I\) satisfies \(\widehat{\Psi}_{2\delta}=\Psi_\delta\circ\Psi_\delta\). Clearly \(\widehat{\Psi}_{2\delta}\) has the required \(Id+(2\delta) \widehat\Delta\) structure, where \(\widehat\Delta\) is differentiable, due to smoothness of \(f\), and therefore locally Lispchitz. The consistency of \eqref{eq:alg1}--\eqref{eq:alg3} implies  that the new integrator is also consistent (see \cite[Theorem II.4.1]{geom}).
In addition, the time-reversibility in Proposition \ref{prop:alg} implies that its order of consistency has to be \emph{even} (see e.g.\ \cite[Section II.3, Section V.1]{geom}); simple examples (say \(f(x) = \|x\|^2\)) show that the order is exactly two.
(An aside: \(\Psi_\delta\) is consistent and time-reversible and yet is not consistent of even order. There is no contradiction, because the argument that shows even order for reversible integrators in \cite[Section II.3, Section V.1]{geom} requires \(\Psi_\delta= Id+\mathcal{O}(\delta)\) something that does not hold for
\eqref{eq:alg1}--\eqref{eq:alg3}.)

For  \(\widehat{\mathcal I}\),
consistency of order two  entails convergence of order two following the standard argument. The even numbered points \((x_{2k},v_{2k})\) generated by \eqref{eq:alg1}--
\eqref{eq:alg3} may be seen as coming from \(\widehat{\mathcal I}\) and this proves
\[
\max_{2k\delta\leq T} \Big(\|x_{2k}-x(2k\delta)\|+\|v_{2k}-v(2k\delta)\|\Big) = \mathcal{O}(\delta^2),\qquad \delta\rightarrow 0.
\]

For an odd-numbered step point \((2k+1)\delta\), one may see the timestepping with  \eqref{eq:alg1}--\eqref{eq:alg3} as being performed into two parts. First, one timesteps from \(t=0\) to
\(t=2k\delta\) and then one performs a single timestep from \(t=2k\delta\) to \(t=(2k+1)\delta\). The first part, as  shown in the last display, introduces \(\mathcal{O}(\delta^2)\) errors. The second part just introduces a single local error, which we know from Theorem \ref{th:consistency} is of size \(\mathcal{O}(\delta^2)\). This concludes the proof.

It may be of some interest to emphasize that since \eqref{eq:alg1}--\eqref{eq:alg3} is consistent of order one and \(\widehat{\mathcal I}\) is consistent of order two, when timestepping with \eqref{eq:alg1}--\eqref{eq:alg3}, the truncation error at the present step will almost be cancelled by the truncation error at the next step. (This is \emph{not} the zig-zagging of the midway points \(x_{k+1/2}\) in Figure~\ref{fig:hug}---that zig-zagging cancels deviations of \(x_{k+1/2}\) and \(x_{k+3/2}\) normal to the level set of \(f\).) The cancellation of truncation errors leading to supraconvergence is borne out in Table~\ref{table:cancellation}, where the columns give the truncation errors in \(x\) at the initial point \((x_0,v_0)\) for \eqref{eq:alg1}--\eqref{eq:alg3} and for \(\widehat{\mathcal I}\).

A similar phenomenon takes place implicitly in Property 6.\ of Proposition~\ref{prop:alg}. In \eqref{eq:fincrement} we saw that the change in \(f\) over a single timestep is \(\mathcal{O}(\delta^2)\) (which matches the first-order consistency of  \eqref{eq:alg1}--\eqref{eq:alg3}). This may wrongly suggest that after \(K=\mathcal{O}(1/\delta)\) steps the change in \(f\) would be \(\mathcal{O}(\delta)\), rather than  \(\mathcal{O}(\delta^2)\) as proved in Proposition \ref{prop:alg}.
 Note that in \eqref{eq:fincrement}, while \(H_{k+1}^-\) and \(H_k^+\) differ by a \(\mathcal{O}(\delta)\) amount, the same is not true of \(v_k\) and \(v_{k+1}\), which differ by a \(\mathcal{O}(1)\) amount, so that \(f_{k+1}-f_k =\mathcal{O}(\delta^2)\) but \(f_{k+1}-f_k \neq \mathcal{O}(\delta^3)\). However, over two consecutive steps
\begin{eqnarray*}
f_{k+2}-f_k &=& \frac{\delta^2}{4}  \Big( H_{k+2}^-[v_{k+2},v_{k+2}]-H_{[k+1}^+[v_{k+1},v_{k+1}]\Big)\\
&& \qquad\qquad +\frac{\delta^2}{4}  \Big( H_{k+1}^-[v_{k+1},v_{k+1}]-H_k^+[v_k,v_k]\Big)\\
&=& \frac{\delta^2}{4}\Big( H_{k+2}^-[v_{k+2},v_{k+2}]-H_k^+[v_k,v_k]\Big)\\
&&\qquad\qquad +\frac{\delta^2}{4}(H_{k+1}^- -H_{k+1}^+)[v_{k+1},v_{k+1}].
\end{eqnarray*}
Here  \(v_{k+2}\) and \( v_k\) differ an by a \(\mathcal{O}(\delta)\) amount, and the same is true for \( H_{k+2}^- \) and \(H_k^+\) and for  \(H_{k+1}^-\)  and \(H_{k+1}^+\). Therefore \(f_{k+2}-f_k = \mathcal{O}(\delta^3)\) (which of course matches the fact that \(\widehat{\mathcal I}\) is consistent of order two).

\bigskip

{\bf Acknowledgements.}  CA was supported in part by a Simons Fellowship while visiting the Newton Institute programme ``Stochastic systems for anomalous diffusion''. JMS has been funded by Ministerio de Ciencia e Innovaci\'{o}n (Spain), project PID2022-136585NB-C21, MCIN/AEI/10.13039/501100011033/FED\-ER, UE.
The authors are very thankful to the  Heilbronn Institute for its support.
\bibliographystyle{abbrv}
\bibliography{bib-hug}

\begin{thebibliography}{10}

\bibitem{beskos-pillai-2013}
A.~Beskos, N.~Pillai, G.~Roberts, J.-M. Sanz-Serna, and A.~Stuart.
\newblock {Optimal tuning of the hybrid Monte Carlo algorithm}.
\newblock {\em Bernoulli}, 19(5A):1501 -- 1534, 2013.

\bibitem{blanes2013}
S.~Blanes, F.~Casas, A.~Farres, J.~Laskar, J.~Makazaga, and A.~Murua.
\newblock New families of symplectic splitting methods for numerical
  integration in dynamical astronomy.
\newblock {\em Applied Numerical Mathematics}, 66:58--72, 2013.

\bibitem{bouchard2018bouncy}
A.~Bouchard-C{\^o}t{\'e}, S.~J. Vollmer, and A.~Doucet.
\newblock The bouncy particle sampler: A nonreversible rejection-free {M}arkov
  chain {M}onte {C}arlo method.
\newblock {\em Journal of the American Statistical Association},
  113(522):855--867, 2018.

\bibitem{duane1987hybrid}
S.~Duane, A.~D. Kennedy, B.~J. Pendleton, and D.~Roweth.
\newblock Hybrid {M}onte {C}arlo.
\newblock {\em Physics Letters B}, 195(2):216--222, 1987.

\bibitem{C24}
M.~C. Escudero.
\newblock {\em Approximate Manifold Sampling}.
\newblock PhD thesis, School of Mathematics, University of Bristol, February
  2024.

\bibitem{youhan-serna-2014}
Y.~Fang, J.~M. Sanz-Serna, and R.~D. Skeel.
\newblock {Compressible generalized hybrid Monte Carlo}.
\newblock {\em The Journal of Chemical Physics}, 140(17):174108, 05 2014.

\bibitem{GS91}
B.~Garcia-Archilla and J.~Sanz-Serna.
\newblock A finite difference formula for the discretization of $d^3/dx^3$ on
  nonuniform grids.
\newblock {\em Mathematics of Computation}, 57(195):239--257, 1991.

\bibitem{GP}
G.~H. Golub and V.~Pereyra.
\newblock The differentiation of pseudo-inverses and nonlinear least squares
  problems whose variables separate.
\newblock {\em SIAM Journal on Numerical Analysis}, 10(2):413--432, 1973.

\bibitem{geom}
E.~Hairer, C.~Lubich, and G.~Wanner.
\newblock {\em Geometric {N}umerical {I}ntegration}.
\newblock Springer Series in Computational Mathematics. Springer-Verlag,
  Berlin, second edition, 2006.

\bibitem{HNW93}
E.~Hairer, S.~N{\o}rsett, and G.~Wanner.
\newblock {\em Solving Ordinary Differential Equations I: Nonstiff Problems}.
\newblock Springer Series in Computational Mathematics. Springer Berlin
  Heidelberg, 1993.

\bibitem{H}
P.~Henrici.
\newblock {\em Discrete {V}ariable {M}ethods in {O}rdinary {D}ifferential
  {E}quations}.
\newblock New York: Wiley, 1962.

\bibitem{Leimkuhler_Reich_2005}
B.~Leimkuhler and S.~Reich.
\newblock {\em Simulating Hamiltonian Dynamics}.
\newblock Cambridge Monographs on Applied and Computational Mathematics.
  Cambridge University Press, 2005.

\bibitem{LS23}
M.~Ludkin and C.~Sherlock.
\newblock Hug and hop: a discrete-time, nonreversible {M}arkov chain {M}onte
  {C}arlo algorithm.
\newblock {\em Biometrika}, 110(2):301--318, 2023.

\bibitem{neal2003slice}
R.~M. Neal.
\newblock Slice sampling.
\newblock {\em The annals of statistics}, 31(3):705--767, 2003.

\bibitem{peters2012rejection}
E.~A. Peters and G.~de~With.
\newblock Rejection-free {M}onte {C}arlo sampling for general potentials.
\newblock {\em Physical Review E---Statistical, Nonlinear, and Soft Matter
  Physics}, 85(2):026703, 2012.

\bibitem{sherlock2021discretebouncyparticlesampler}
C.~Sherlock and A.~H. Thiery.
\newblock A discrete bouncy particle sampler.
\newblock {\em arXiv preprint arXiv:1707.05200}, 2017.

\bibitem{sherlock2022discrete}
C.~Sherlock and A.~H. Thiery.
\newblock A discrete bouncy particle sampler.
\newblock {\em Biometrika}, 109(2):335--349, 2022.

\bibitem{vanetti2019piecewise}
P.~Vanetti.
\newblock {\em Piecewise-deterministic Markov Chain Monte Carlo}.
\newblock PhD thesis, University of Oxford, 2019.

\bibitem{vanetti2017piecewise}
P.~Vanetti, A.~Bouchard-C{\^o}t{\'e}, G.~Deligiannidis, and A.~Doucet.
\newblock Piecewise-deterministic {M}arkov chain {M}onte {C}arlo.
\newblock {\em arXiv preprint arXiv:1707.05296}, 2017.

\bibitem{zappa-holmes}
E.~Zappa, M.~Holmes-Cerfon, and J.~Goodman.
\newblock Monte {C}arlo on manifolds: Sampling densities and integrating
  functions.
\newblock {\em Communications on Pure and Applied Mathematics}, 71, 02 2017.

\end{thebibliography}



\end{document}